# RANDOMLY GROWING BRAID ON THREE STRANDS AND THE MANTA RAY

By Jean Mairesse and Frédéric Mathéus

*CNRS–Université Paris 7 and Université de Bretagne-Sud*

*In memory of Daniel Mollier, our former mathematics teacher at the Lycée Louis le Grand, Paris.*

Consider the braid group $B_3 = \langle a, b | aba = bab \rangle$ and the nearest neighbor random walk defined by a probability $\nu$ with support $\{a, a^{-1}, b, b^{-1}\}$. The rate of escape of the walk is explicitly expressed in function of the unique solution of a set of eight polynomial equations of degree three over eight indeterminates. We also explicitly describe the harmonic measure of the induced random walk on $B_3$ quotiented by its center. The method and results apply, mutatis mutandis, to nearest neighbor random walks on dihedral Artin groups.

**1. Introduction.** Consider three strands and the following four elementary moves: $a =$ cross the strands 1 and 2 with 1 above; $a^{-1} =$ cross the strands 1 and 2 with 2 above; $b =$ cross the strands 2 and 3 with 2 above; $b^{-1} =$ cross the strands 2 and 3 with 3 above.

At any given instant, a strand is numbered from 1 to 3 according to the position of its bottom extremity from left to right. Starting from an uncrossed configuration and applying elementary moves successively, we obtain a *braid diagram* as in Figure 1.

Two braid diagrams are equivalent if they are projections of isotopic three-dimensional configurations, that is, if we can go from one diagram to the other by moving around the strands without touching the top and bottom extremities. For instance, the three diagrams of Figure 1 are equivalent. A *braid* is an equivalence class of braid diagrams. The *length* $|X|$ of a braid $X$ is the minimal number of elementary moves necessary to build the braid. For instance, the length of the braid of Figure 1 is 3.

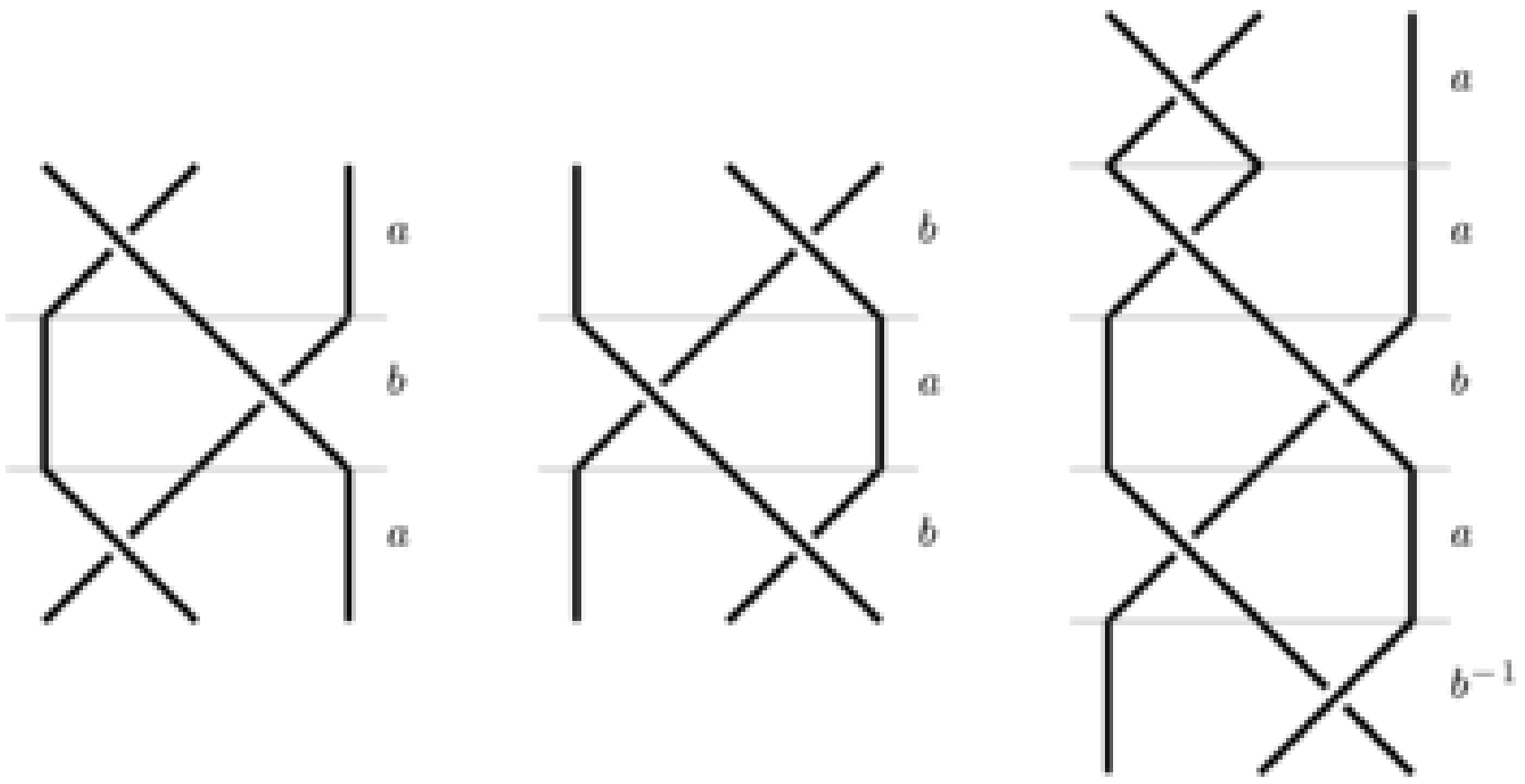

FIG. 1. *From left to right, the braid diagrams obtained by performing the moves* $(a, b, a)$, $(b, a, b)$ *and* $(a, a, b, a, b^{-1})$.

A *randomly growing braid* is the random process $(X_n)_n$ obtained by choosing the successive elementary moves independently and according to some given probability distribution $\nu$ on $\{a, a^{-1}, b, b^{-1}\}$. It has been proposed as a pertinent model to study the random entanglement of polymers or chain-like objects [26, 27]. Other applications in fluid mechanics and astrophysics are described in [2]. A main quantity of interest is the growth rate of the length: $\gamma = \lim_n |X_n|/n$. It captures the speed at which the strands get entangled together.

The above problem can be restated in more algebraic terms. Given two braids $u$ and $v$, the braid $uv$ is obtained by "gluing" a braid diagram of $u$ on top of a braid diagram of $v$. The set of braids with the above internal law forms a group which can be proved to have the finite group presentation

$$B_3 = \langle a, b | aba = bab \rangle.$$

The group $B_3$ is known as the *braid group on three strands*, and the generators $\{a, a^{-1}, b, b^{-1}\}$ as the Artin generators. The random process $(X_n)_n$ defined above is a realization of the nearest neighbor random walk (NNRW) on the group $B_3$ defined by the probability $\nu$ on the generators $\{a, a^{-1}, b, b^{-1}\}$. The computed limit $\gamma$ is the drift or rate of escape of the walk. The way to determine the drift is by expliciting a higher level object, which is the harmonic measure (or exit measure) of the random walk. Roughly, it gives the direction taken by the walk in its escape to infinity. The harmonic measure is important since it captures a lot of information about the RW.

In this paper we compute explicitly $\gamma$ for any distribution $\nu$ on $\{a, a^{-1}, b, b^{-1}\}$. More precisely, given $\nu$, we define eight polynomial equations of degree 3 over eight indeterminates [see (30)], which are shown to admit a unique solution $r$. The rate $\gamma$ is then obtained as an explicit linear functional of $r$; see (29). The polynomial equations can be partially or completely solved to provide a closed form formula for $\gamma$ under the following two natural symmetries: (i) $\nu(a)=\nu(a^{-1})$ and $\nu(b)=\nu(b^{-1})$; (ii) $\nu(a)=\nu(b)$ and $\nu(a^{-1})=\nu(b^{-1})$. In case (i), for $p=\nu(a)=\nu(a^{-1})\in(0,1/2)$,

$$\gamma = p+(1-4p)u,$$

with

$$2(4p-1)u^3+(24p^2-18p+1)u^2+p(7-12p)u+p(2p-1)=0,$$

and $u$ is the smallest positive root of the polynomial. In case (ii), for $p=\nu(a)=\nu(b)\in(0,1/2)$,

$$\gamma=\max\biggl[1-4p,\frac{(1-2p)(-1-4p+\sqrt{5-8p+16p^2})}{2(1-4p)},\frac{p(-3+4p+\sqrt{5-8p+16p^2})}{-1+4p},-1+4p\biggr]. \quad (1)$$

Assume that $p=\nu(a)=1/2-\nu(b^{-1})$ and $q=\nu(b)=1/2-\nu(a^{-1})$. This is satisfied both in case (i) and in case (ii). Under these restrictions, we can have a three-dimensional representation of $\gamma$ as a function of $p$ and $q$. The result, obtained with the help of Maple, is given in Figure 2.

Since a set of algebraic equations can be solved with arbitrary precision, Figure 2 should be considered as an exact picture rather than an approximated one. Observe that the surface exhibits nondifferentiability lines. They correspond to some phase transitions whose physical meaning is intriguing. Let us give some partial explanation of the phenomenon.

Consider the surface of Figure 2 from above. The nondifferentiability lines partition the square into four regions, as shown on Figure 3 (left). In the corner regions, the limiting random braid contains only the positive generators $a$ and $b$ (right-upper corner), or only the negative generators $a^{-1}$ and $b^{-1}$ (left-lower corner). It implies that the drift is simply $\gamma=|\nu(a)+\nu(b)-\nu(a^{-1})-\nu(b^{-1})|=|2p+2q-1|$. In the other two regions, there is a coexistence of positive and negative generators in the limit, with more positive than negative in the right region, and the other way round in the left region.

Another perspective on these four regions is obtained by comparing the RW on $B_3$ with its projection on $\mathbb{Z}$. Assume, for simplicity, that $p=q$ [case (ii)]. Let $\varphi : B_3\to\mathbb{Z}$ be the group morphism defined by $\varphi(a)=\varphi(b)=1$, and $\varphi(a^{-1})=\varphi(b^{-1})=-1$. Then $(\varphi(X_n))_n$ is a realization of the NNRW on $\mathbb{Z}$ defined by the probability $\lambda : \lambda(1)=\nu(a)+\nu(b)=2p$, $\lambda(-1)=1-2p$. So the drift of $(\varphi(X_n))_n$ is $\gamma_{\mathbb{Z}}=|4p-1|$, with a point of nondifferentiability at

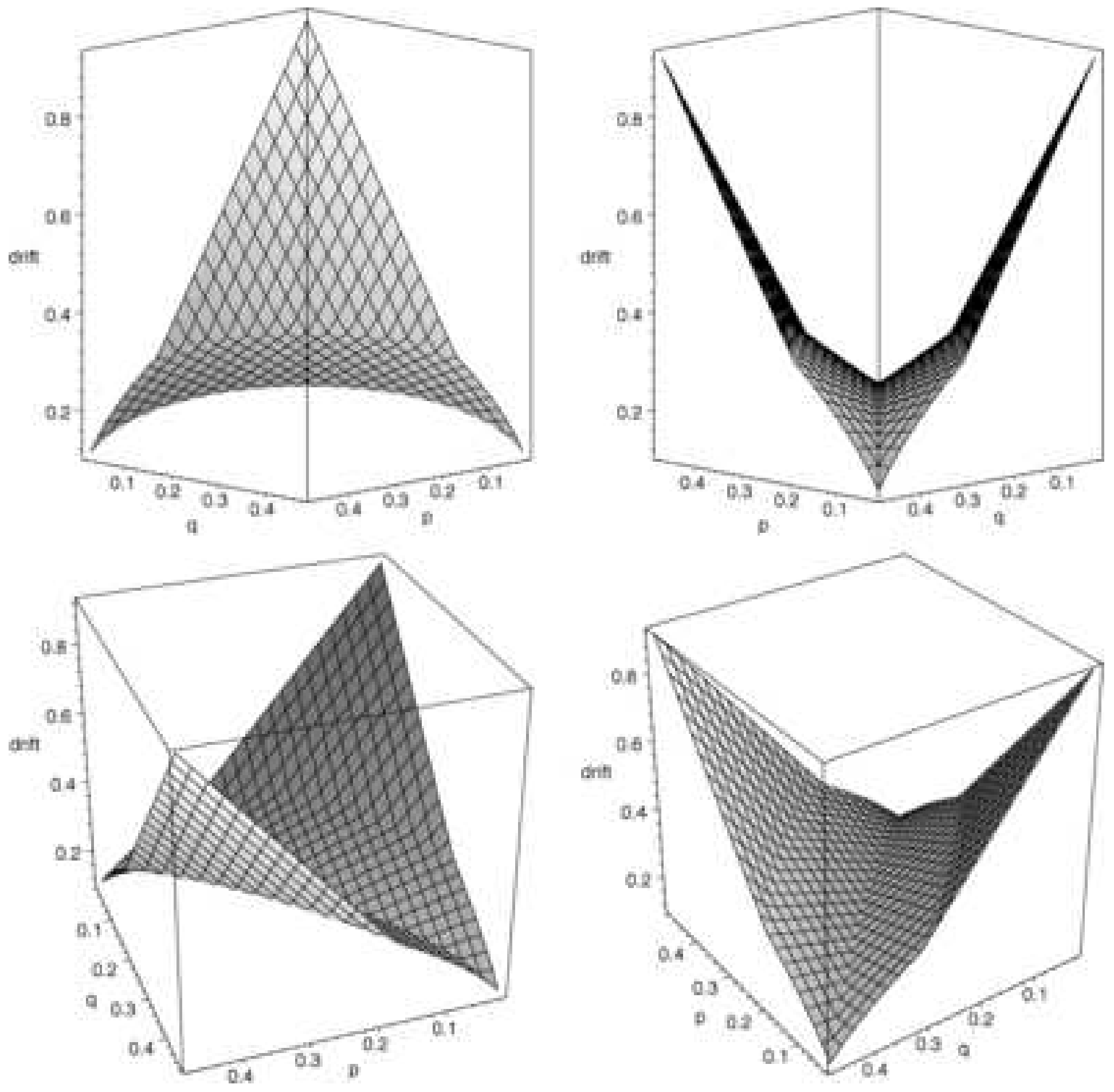

FIG. 2. *The manta ray from four different angles. That is, the drift* $\gamma$ *as a function of* $p = \nu(a) = 1/2 - \nu(b^{-1})$ *and* $q = \nu(b) = 1/2 - \nu(a^{-1})$.

$p = 1/4$. The interpretation of the corresponding phase transition is easy: the random walker escapes to infinity along the positive axis for $2p > 1/2$, and along the negative axis for $2p < 1/2$. No such simple explanation is available for $B_3$. However, the three phase transitions for $B_3$ are related to the projected RW on $\mathbb{Z}$. In Figure 3 (right) we have superposed the drift curves for $B_3$ ($\gamma_{B_3}$ in black) and for $\mathbb{Z}$ ($\gamma_{\mathbb{Z}}$ in light gray). The points of nondifferentiability of $\gamma_{B_3}$ are $p = 1/4$, the same as $\gamma_{\mathbb{Z}}$, and also the two points where the curves $\gamma_{B_3}$ and $\gamma_{\mathbb{Z}}$ separate. One can observe [or check numerically using (1)] that even the difference $\gamma_{B_3} - \gamma_{\mathbb{Z}}$ is nondifferentiable at $p = 1/4$, which may also be surprising.

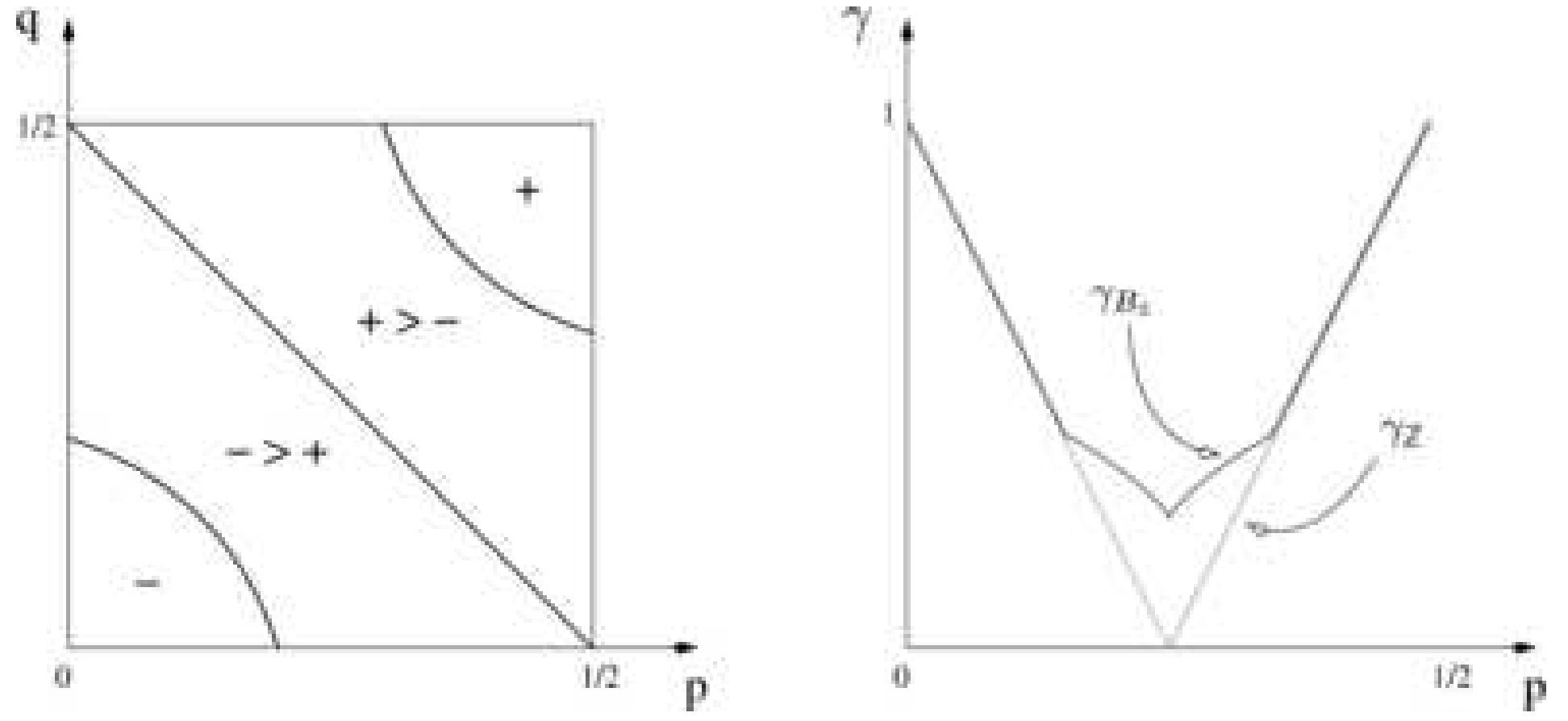

FIG. 3. *Phase transitions for the random braid model.*

Combinatorial and probabilistic arguments are entangled to get the above results. The central idea is that the randomly growing Garside normal form of a braid is well behaved. In a nutshell, the different steps of the proofs are as follows:

(A) Go from $B_3$ to $B_3/Z$, the quotient by the center. Switch from the Artin generators $S$ to the Garside generators $T$.

(B) Consider the induced random walk on $(B_3/Z,T)$. Show that the harmonic measure $\mu^\infty$ can be explicitly computed. It has a multiplicative structure of the type

$$\mu^\infty(u_1u_2\cdots u_kT^{\mathbb{N}}) = \alpha\mathcal{M}(u_1)\cdots\mathcal{M}(u_{k-1})\beta(u_k),$$

with $\alpha \in \mathbb{R}_+^{1\times 4}, \mathcal{M}:T\to\mathbb{R}_+^{4\times 4}, \beta:T\to\mathbb{R}_+^{4\times 1}$. The coefficients of $\alpha,\mathcal{M}(\cdot),\beta$ have a simple expression in terms of the unique solution to a set of polynomial equations. Based on this, compute the drift on $B_3/Z$ as well as the drift along the center.

(C) Go back from $B_3/Z$ to $B_3$. Switch back from $T$ to $S$. This very last step requires combinatorial work since the geodesics according to $T$ are quite different from the geodesics according to $S$.

Let us insist on one point concerning step B. It is well known that the group $B_3/Z$ is isomorphic to the modular group $\mathbb{Z}/2\mathbb{Z}\star\mathbb{Z}/3\mathbb{Z}$. Precisely, we have $B_3/Z \sim \langle u,v|u^2=v^3=1\rangle$ with $u$ corresponding to $aba$ and $v$ to $ab$. NNRW on $\mathbb{Z}/2\mathbb{Z}\star\mathbb{Z}/3\mathbb{Z}$ are specifically studied in [21], Section 4.2, and the harmonic measure and the drift are given in closed form. Here the model is different since we need to deal with the set of generators $S$ or $T$. The Cayley graph of $(B_3/Z,T)$ is more twisted than the one of $(\mathbb{Z}/2\mathbb{Z}\star\mathbb{Z}/3\mathbb{Z},\{u,v,v^{-1}\})$ (see Figures 4 and 5), and the NNRW is consequently more complex to analyze.

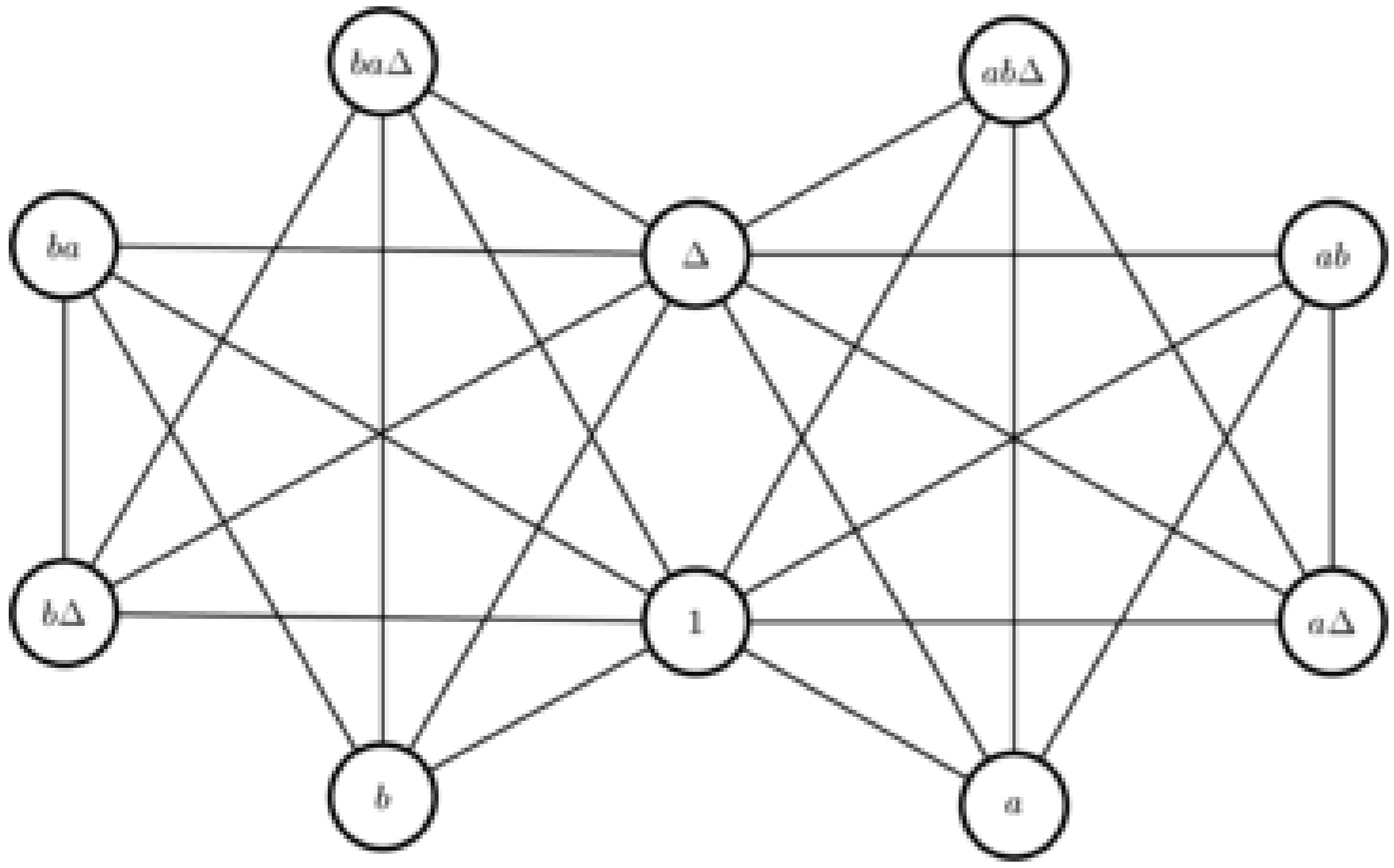

Fig. 4. *The Cayley graph* $\mathcal{X}(B_3/Z, \Sigma)$ *restricted to the nodes* $\Sigma \cup \{1, \Delta\}$.

There are two natural families of groups generalizing $B_3$: the braid groups $B_n$ and the Artin groups of dihedral type $I_2(n)$. For the braid groups $B_n, n \geq 4$, the above approach completely collapses. First, it is not clear that the random Garside normal form should have good asymptotic properties. Second, switching back from Garside to Artin generators is a priori out of reach. Understanding the combinatorics of Artin geodesics is indeed a well-known open and difficult question. It is not even known if the growth series of $B_n$ with respect to the Artin generators is rational. On the other hand, the results and method of proof for $B_3$ extend to Artin groups of dihedral type. We treat this extension in Section 5.

The approach can also be extended to a randomly growing braid on 3 strands on a "cylinder" [see Section 4.4] or to the dual structure of dihedral Artin groups.

Random walks on discrete infinite groups form a vast and active domain; see the monograph by Woess [31], or for a different perspective, the forthcoming monograph by Lyons and Peres [18]. In this landscape, the contribution of the paper is as follows. We provide a complete and detailed study of the drift for a specific group, $B_3$, which is not virtually free and whose combinatorics is nonelementary. In the past not too much effort has been devoted to getting such explicit computations. Here is a presumably nonexhaustive list of papers where nonelementary computations are carried out: [6, 7, 15, 17, 19, 21, 25, 27, 29, 30]. The articles [27, 30] are of particular relevance in our context since they deal with braid groups.

There are also two general potential methods for computing the drift. The first one is due to Sawyer and Steger [29], it was developped for NNRW on homogeneous trees but adapts to NNRW on plain groups (i.e., free products of free and finite groups). In this approach the drift is expressed as a functional of the *first-passage generating series* of the random walk. (Concerning the computation of the generating series of the walk, there exists an important literature; see [31], Sections II.9 and III.17.) The second method was proposed in [19, 21], and also applies to NNRW on plain groups. There, the drift is expressed as a function of the one-dimensional marginals of the harmonic measure which can be simply derived. Here we proceed along the lines of the second approach, but an important refinement and adaptation of the method is required.

Propositions 4.5 and 5.4 were announced without proofs in [20]. Part of the computations, as well as a couple of extensions, are not detailed in the present paper. They can be found in the version with the Appendix posted on the arXiv [23].

**2. Preliminaries.** Let $\mathbb{N}$ be the set of nonnegative integers and let $\mathbb{N}^* = \mathbb{N}\backslash\{0\}$. We also set $\mathbb{R}_+^* = \mathbb{R}_+\backslash\{0\}$. We denote the support of a random variable by supp. If $\mu$ is a measure on a group $(G, *)$, then $\mu^{*n}$ is the $n$-fold convolution product of $\mu$, that is, the image of the product measure $\mu^{\otimes n}$ by the product map $G \times \cdots \times G \to G, (g_1, \dots, g_n) \mapsto g_1 * g_2 * \cdots * g_n$.

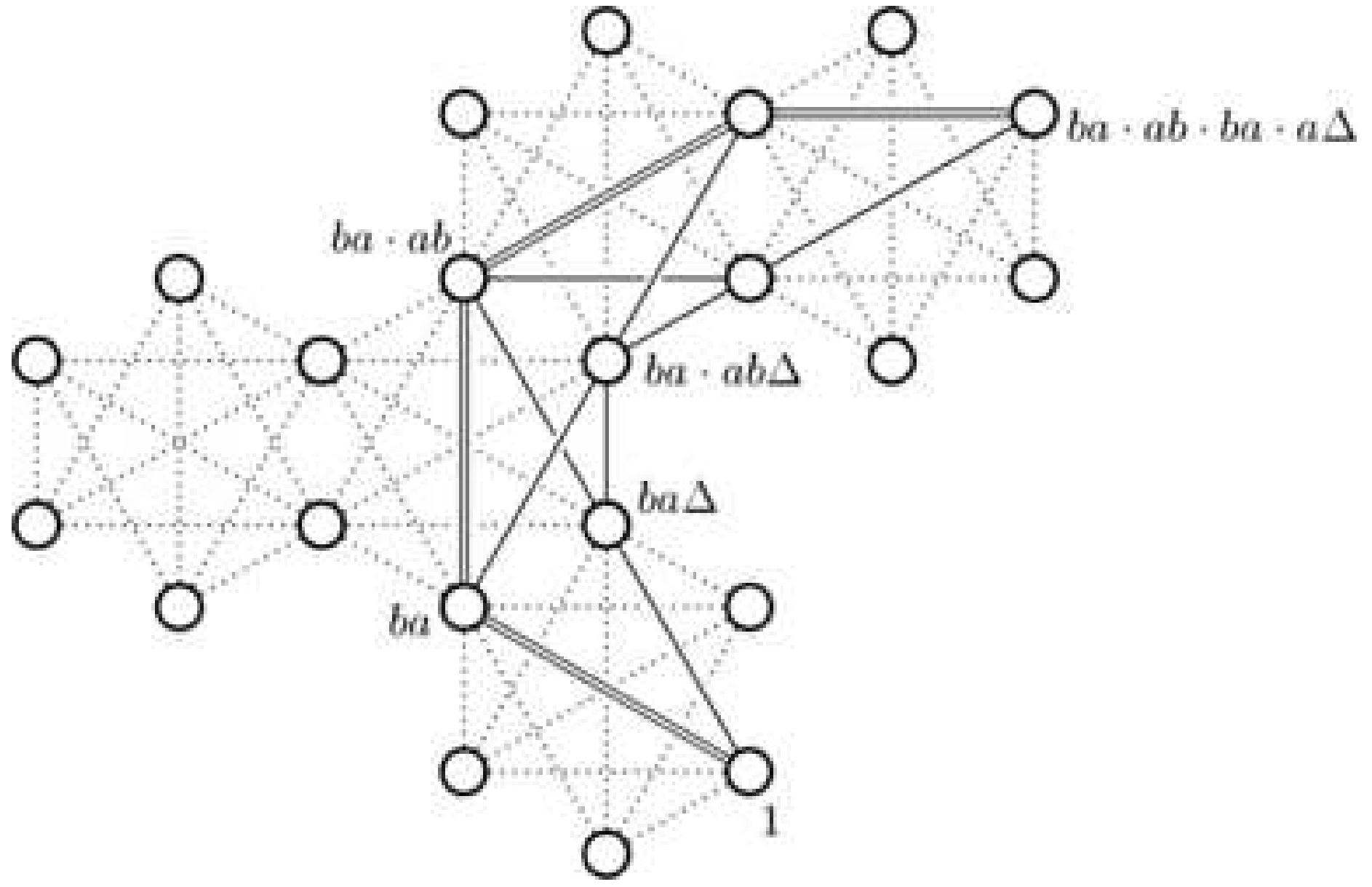

FIG. 5. *The Cayley graph $\mathcal{X}(B_3/Z, \Sigma)$.*

The symbol $\sqcup$ is used for the disjoint union of sets. Given a finite set $\Sigma$, a vector $x \in \mathbb{R}^\Sigma$ and $S \subset \Sigma$, set $x(S) = \sum_{u\in S} x(u)$. The power set of a set $S$ is denoted by $\mathcal{P}(S)$.

Given a set $\Sigma$, the free monoid it generates is denoted by $\Sigma^*$. The empty word is denoted by 1. The concatenation product of $u$ and $v$ is written $u \cdot v$ or simply $uv$. The *length* (number of letters) of a word $u$ is denoted by $|u|_\Sigma$.

Let $\mathbb{K}$ be a semiring. Let $Q$ and $\Sigma$ be two finite sets. A $\mathbb{K}$-*automaton* of set of states (dimension) $Q$ over the alphabet $\Sigma$ is a triple $\mathcal{A} = (\alpha, \mu, \beta)$, where $\alpha \in \mathbb{K}^{1\times Q}$, $\beta \in \mathbb{K}^{Q\times 1}$, and where $\mu : \Sigma^* \to \mathbb{K}^{Q\times Q}$ is a morphism of monoids. The morphism $\mu$ is uniquely determined by the family of matrices $\{\mu(a), a \in \Sigma\}$, and for $w = a_1 \cdots a_n$, we have $\mu(w) = \mu(a_1)\mu(a_2)\cdots\mu(a_n)$. The (formal power) series *recognized* by $\mathcal{A}$ is $S : \Sigma^* \to \mathbb{K}, S(w) = \alpha\mu(w)\beta$. The set of series recognized by a $\mathbb{K}$-automaton is denoted by $\mathbb{K}\,\mathrm{Rat}(\Sigma^*)$.

An *automaton* is a $\mathbb{B}$-automaton where $\mathbb{B}$ is the Boolean semiring ({false, true}, or, and) $\sim (\{0,1\}, \max, \min)$. In this case, we identify the series $S$ with the language $\{w \in \Sigma^* | S(w) = 1\}$, and we call it the language *recognized* by the automaton. For a finite set $U$, the set $\mathbb{B}\,\mathrm{Rat}(U^*)$, with union and concatenation as the two laws, forms a semiring. A $\mathbb{B}\,\mathrm{Rat}(U^*)$-automaton is called a *transducer*, and the associated series is a *transduction*, that is, a map from $\Sigma^*$ to $\mathcal{P}(U^*)$.

A state $i \in Q$ is *initial*, respectively, *final*, if $\alpha_i \neq \mathbb{0}$, respectively, $\beta_i \neq \mathbb{0}$ (where $\mathbb{0}$ is the zero element of the semiring). As usual, a $\mathbb{K}$-automaton is represented graphically by a labeled and weighted directed graph with ingoing and outgoing arcs for initial and final states. For $\alpha_i = x \neq \mathbb{0}$, $\mu(a)_{jk} = y \neq \mathbb{0}$ and $\beta_l = z \neq \mathbb{0}$, we have the following respective ingoing arc, arc and outgoing arc:

$$\xrightarrow{|x} i, \qquad j \xrightarrow{a|y} k, \qquad l \xrightarrow{|z} .$$

We call $a$ the *label* and $y$ (or $x$, or $z$) the *weight* of the arc. Examples of $\mathbb{K}$-automata appear in Figure 6 ($\mathbb{K} = \mathbb{B}$), Figure 8 ($\mathbb{K} = \mathbb{B}\,\mathrm{Rat}(\Sigma^*)$) or Figure 9 ($\mathbb{K} = (\mathbb{R}_+, +, \times)$). (In the figures, the weights equal to $\mathbb{1}$, the unit element of the semiring, are omitted.) The *label* of a path is the concatenation of the labels of the arcs. A path which is both starting with an ingoing arc and ending with an outgoing arc is called a *successful path*.

For details on automata with multiplicities over a semiring see, for instance, [4, 8].

Consider a group $(G, *)$. The unit element is denoted by 1, and the inverse of the element $u$ is denoted by $u^{-1}$. We often omit the $*$ sign, writing $xy$ instead of $x * y$. Assume $G$ to be finitely generated and let $\Sigma \subset G$ be a finite set of generators of $G$. Denote by $\pi : \Sigma^* \to G$ the monoid homomorphism which associates to a word $a_1 \cdot \cdots \cdot a_k$ the group element $a_1 * \cdots * a_k$. A word $u \in \pi^{-1}(g)$ is called a *representative* of $g$. A language $L$ of $\Sigma^*$ is a *cross-section* or *set of normal forms* of $G$ (over the alphabet $\Sigma$) if the restriction of $\pi$ to $L$ defines a bijection, that is, if every element of $G$ has a unique representative in $L$.

The *length with respect to* $\Sigma$ of a group element $g$ is

$$|g|_\Sigma = \min\{k | g = s_1 * \cdots * s_k, s_i \in \Sigma\}. \tag{2}$$

A representative $u$ of $g$ is a *geodesic* (*word*) if $|u|_\Sigma = |g|_\Sigma$. A *geodesic cross-section* is defined accordingly.

The *Cayley graph* $\mathcal{X}(G,\Sigma)$ of the group $G$ with respect to the set of generators $\Sigma$ is the directed graph with set of nodes $G$ and with an arc from $u$ to $v$ if $\exists a \in \Sigma, u * a = v$.

Let $\mu$ be a probability distribution over $\Sigma$. Consider the Markov chain on the state space $G$ with one-step transition probabilities given by $\forall g \in G, \forall a \in \Sigma, P_{g,g*a} = \mu(a)$. This Markov chain is called the *random walk* (associated with) $(G,\mu)$. It is a *nearest neighbor* random walk: one-step moves occur between nearest neighbors in the Cayley graph $\mathcal{X}(G,\Sigma)$. Let $(x_n)_n$ be a sequence of i.i.d. r.v.'s distributed according to $\mu$. Set

$$X_0 = 1, \qquad X_{n+1} = X_n * x_n = x_0 * x_1 * \cdots * x_n. \tag{3}$$

The sequence $(X_n)_n$ is a *realization* of the random walk $(G,\mu)$. The law of $X_n$ is $\mu^{*n}$. By subadditivity [12], there exists a constant $\gamma \in \mathbb{R}_+$ such that a.s. and in $L^p$, for all $1 \le p < \infty$,

$$\lim_{n\to\infty} \frac{|X_n|_\Sigma}{n} = \gamma. \tag{4}$$

We call $\gamma$ the *drift* or *rate of escape* of the walk.

If the group $G$ is nonamenable and if $\mu$ generates the whole group, that is, $\bigcup_{n\in\mathbb{N}} \operatorname{supp}(\mu^{*n}) = G$, then the random walk $(G,\mu)$ is transient and the drift is positive (see [12] and [31], Chapter 1.B for details). Below, all the groups considered are nonamenable so we only deal with transient random walks.

**3. The braid group on three strands.** The *braid group on three strands* is the finitely-presented group $B_3 = \langle a, b | aba = bab\rangle$. The set $S = \{a, a^{-1}, b, b^{-1}\}$ is the set of *Artin* generators. The name of the group comes from the usual graphical interpretation in terms of braids; see Figure 1.

Set $\Delta = aba = bab$. Observe that $a\Delta = a(bab) = (aba)b = \Delta b$. Similarly, we have $b\Delta = \Delta a, a^{-1}\Delta = \Delta b^{-1}, b^{-1}\Delta = \Delta a^{-1}$. The center of $B_3$ is precisely the subgroup:

$$Z = \langle \Delta^2 \rangle = \{\Delta^{2k}, k \in \mathbb{Z}\}.$$

Denote by $B_3/Z$ the quotient of $B_3$ by $Z$. Define the canonical homomorphism $p : B_3 \to B_3/Z$. Given an element $u$ in $B_3$, its class $p(u)$ in $B_3/Z$ is still denoted by $u$, with the interpretation depending on the context. We have the following:

PROPOSITION 3.1 (Garside [10]). *Set* $T = \{a, b, ab, ba\}$. *For* $g \in B_3$, *there exist* $g_1, \ldots, g_m \in T$ *and* $k \in \mathbb{Z}$ *such that*

$$g = g_1 g_2 \cdots g_m \Delta^k.$$

*If* $m$ *is chosen to be minimal, then the decomposition is unique.*

Given $g = g_1 g_2 \cdots g_m \Delta^k \in B_3$, we have, in $B_3/Z$, $p(g) = g_1 \cdots g_m$ if $k$ is even, and $p(g) = g_1 \cdots g_m \Delta$ if $k$ is odd. So we can restate Proposition 3.1, replacing "$g \in B_3$" by "$g \in B_3/Z$" and "$k \in \mathbb{Z}$" by "$k \in \{0, 1\}$."

Both in $B_3$ and $B_3/Z$, the unique decomposition corresponding to a minimal value of $m$ is called the *Garside normal form* of $g$.

In $B_3$ it is natural to view the Garside normal form as an element of $V^*$ with $V = \{a, b, ab, ba, \Delta, \Delta^{-1}\}$. Given $g \in B_3$ with Garside normal form $g_1 \cdots g_m \Delta^k$, one can prove that $|g|_V = m + |k|$, meaning that the set of Garside normal forms is a geodesic cross-section over $V$. In $B_3/Z$ we adopt a different point of view by considering the a priori less natural set of generators $\Sigma$ defined below. In $B_3/Z$ we have $\Delta^{-1} = \Delta$ and,

$$a^{-1} = ba\Delta, \qquad b^{-1} = ab\Delta, \qquad (ab)^{-1} = a\Delta, \qquad (ba)^{-1} = b\Delta.$$

Define the following subsets of $B_3/Z$:

$$T = \{a, b, ab, ba\}, \qquad \Sigma = T \cup T^{-1} = \{a, b, ab, ba, a\Delta, b\Delta, ab\Delta, ba\Delta\}.$$

Observe that $\Sigma$ is a set of generators of $B_3/Z$. Define the canonical homomorphism:

$$\pi : \Sigma^* \to B_3/Z. \tag{5}$$

In Figures 4 and 5 we have represented the Cayley graph $\mathcal{X}(B_3/Z, \Sigma)$.

Let $\mathcal{G}$ be the set of Garside normal forms of $B_3/Z$. We view $\mathcal{G}$ as being a subset of $\Sigma^*$, and more precisely, of $T^*\Sigma$: an element $g_1 \cdots g_m \Delta^{0/1}, m > 0$, is viewed as $g_1 \cdot \cdots \cdot g_{m-1} \cdot (g_m \Delta^{0/1})$. For $\Delta$, we arbitrarily choose a geodesic representative in $\Sigma^*$, say $a \cdot ba$, that we still denote by $\Delta$ for simplicity. It follows easily from the definition of the Garside normal form that $\mathcal{G}$ is a geodesic cross-section over $\Sigma$ of $B_3/Z$. Let

$$\phi : B_3/Z \to \mathcal{G} \tag{6}$$

be the map which associates to a group element its Garside normal form.

The set $\mathcal{G}$ belongs to $\mathrm{Rat}(\Sigma^*)$. This is an easy consequence of the proof by Thurston [9], Chapter 9, that braid groups are automatic. Let us detail this point.

Define the two maps $\mathrm{First}, \mathrm{Last} : \Sigma \to \{a, b\}$ by Table 1.

Observe that $\mathrm{First}(u)$ is the first symbol of $u$. For $u \in T$, $\mathrm{Last}(u)$ is the last symbol of $u$ and $\mathrm{Last}(u\Delta)$ is the last symbol of $v$, where $\Delta v = u\Delta$.

TABLE 1

| | $a$ | $b$ | $ab$ | $ba$ | $a\Delta$ | $b\Delta$ | $ab\Delta$ | $ba\Delta$ |
|---|---|---|---|---|---|---|---|---|
| First | $a$ | $b$ | $a$ | $b$ | $a$ | $b$ | $a$ | $b$ |
| Last | $a$ | $b$ | $b$ | $a$ | $b$ | $a$ | $a$ | $b$ |

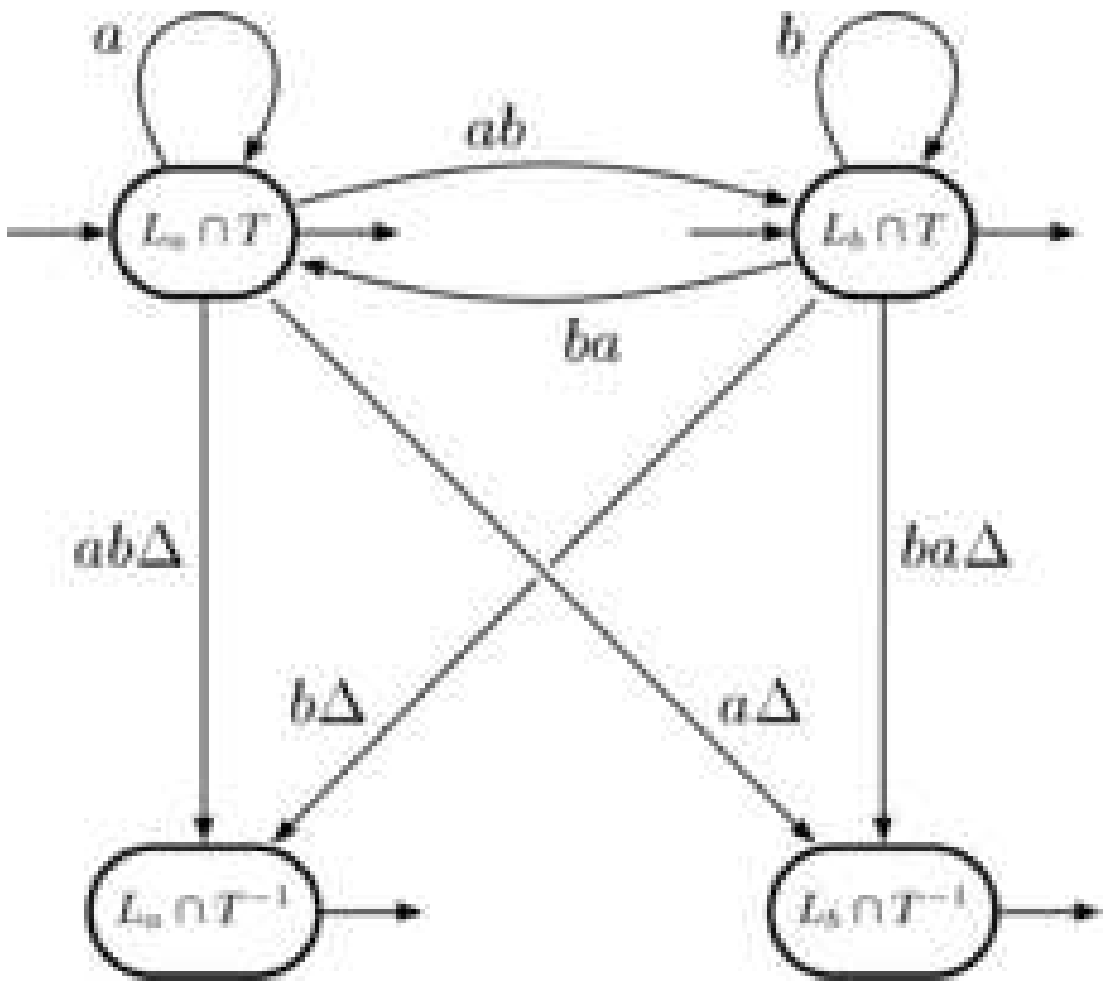

FIG. 6. *Automaton over $\Sigma$ recognizing the Garside geodesics $\mathcal{G}$.*

Set $L_a = \text{Last}^{-1}(a), L_b = \text{Last}^{-1}(b)$. Consider $u_1 \cdots u_k \in \Sigma^*$. We have [9], Chapter 9

$$u_1 \cdots u_k \in \mathcal{G} \backslash \{\Delta\} \quad \Longleftrightarrow \quad u_1 \cdots u_k \in T^* \Sigma \qquad (7)$$
$$\forall i \in \{1, \ldots, k-1\},\ \text{Last}(u_i) = \text{First}(u_{i+1}).$$

The element $\Delta = a \cdot ba \in \mathcal{G}$ is specific in that it does not satisfy (7). For simplicity, in Figures 6–8 we have omitted to represent $\Delta$.

Rephrasing (7), we get that $\mathcal{G}$ is recognized by the automaton given in Figure 6. The states are labeled by disjoint subsets partitioning $\Sigma$. A successful path ending in a state with label $L$ is labeled by an element of $\mathcal{G}$ whose last letter belongs to $L$.

By systematic inspection, one can check the following crucial observation:

$$\forall u, v \in \Sigma \qquad w = u * v \in \Sigma \quad \Longrightarrow \quad \text{First}(w) = \text{First}(u), \qquad (8)$$
$$\text{Last}(w) = \text{Last}(v).$$

It is convenient to introduce a notation for the law of $\mathcal{G}$ induced by the group law of $B_3/Z$. Define

$$\begin{aligned} &\mathcal{G} \times \mathcal{G} \to \mathcal{G}, \\ &(u,v) \mapsto u \circledast v = \phi(\pi(u) * \pi(v)). \end{aligned} \tag{9}$$

By definition, $(\mathcal{G}, \circledast)$ is a group which is isomorphic to $(B_3/Z, *)$, the isomorphism being $\phi : B_3/Z \to \mathcal{G}$. It is instructive to give a recursive definition of the law $\circledast$. Let $\iota : \Sigma \to \Sigma$ be the involution defined by

$$\iota(a) = b, \qquad \iota(ab) = ba, \qquad \iota(a\Delta) = b\Delta, \qquad \iota(ab\Delta) = ba\Delta. \tag{10}$$

We have $\forall u \in \mathcal{G}, 1 \circledast u = u \circledast 1 = u$. We have $\Delta \circledast \Delta = 1$ and $\forall u_1 \cdots u_k \in \mathcal{G} \backslash \{1, \Delta\}$,

$$\begin{aligned} &\Delta \circledast (u_1 \cdots u_k) = \iota(u_1) \cdots \iota(u_{k-1}) \cdot [\iota(u_k)\Delta], \\ &(u_1 \cdots u_k) \circledast \Delta = u_1 \cdots u_{k-1} \cdot [u_k \Delta]. \end{aligned}$$

We have $\forall u = u_1 \cdots u_k \in \mathcal{G} \backslash \{1, \Delta\}, \forall v = v_1 \cdots v_l \in \mathcal{G} \backslash \{1, \Delta\}$,

$$u \circledast v = \begin{cases} (u_1 \cdots u_{k-1}) \circledast (v_2 \cdots v_l), \\ \qquad \text{if } u_k * v_1 = 1, \\ (u_1 \cdots u_{k-1}) \circledast (\iota(v_2) \cdots \iota(v_{l-1}) \cdot [\iota(v_l)\Delta]), \\ \qquad \text{if } u_k * v_1 = \Delta, \\ u_1 \cdots u_{k-1} \cdot w \cdot v_2 \cdots v_l, \\ \qquad \text{if } u_k * v_1 = w \in T, \\ u_1 \cdots u_{k-1} \cdot w \cdot \iota(v_2) \cdots \iota(v_{l-1}) \cdot [\iota(v_l)\Delta], \\ \qquad \text{if } u_k * v_1 = w\Delta \in T^{-1}, \\ u_1 \cdots u_k \cdot v_1 \cdots v_l, \\ \qquad \text{if } u_k \in T, u_k * v_1 \notin \Sigma \cup \{1, \Delta\}, \\ u_1 \cdots u_{k-1} \cdot w_k \cdot \iota(v_1) \cdots \iota(v_{l-1}) \cdot [\iota(v_l)\Delta], \\ \qquad \text{if } u_k = w_k\Delta \in T^{-1}, u_k * v_1 \notin \Sigma \cup \{1, \Delta\}. \end{cases} \tag{11}$$

The first four lines correspond to the case [$\text{Last}(u_k) \neq \text{First}(v_1)$] and the last two lines to the case [$\text{Last}(u_k) = \text{First}(v_1)$]. The expression in line 3 is a direct consequence of (8). Part of the complexity in the above definition comes from the fact that the $\Delta$'s have to be shifted to the right in a Garside normal form.

Recall that $\mathcal{G} \subset \Sigma^*$ is a geodesic cross-section of $B_3/Z$. Let $\mathcal{SG} \subset \Sigma^*$ be the set of all the geodesics of $B_3/Z$. Consider $u_1 \cdots u_k \in \Sigma^*$. We have [compare with (7)]

$$u_1 \cdots u_k \in \mathcal{SG} \backslash \{\Delta\} \iff \forall i \in \{1, \dots, k-1\}, \text{Last}(u_i) = \text{First}(u_{i+1}). \tag{12}$$

Again, the element $\Delta$ is specific. It follows from (12) that $\mathcal{SG}$ belongs to $\text{Rat}(\Sigma^*)$ and that the automaton of Figure 7 recognizes $\mathcal{SG}$.

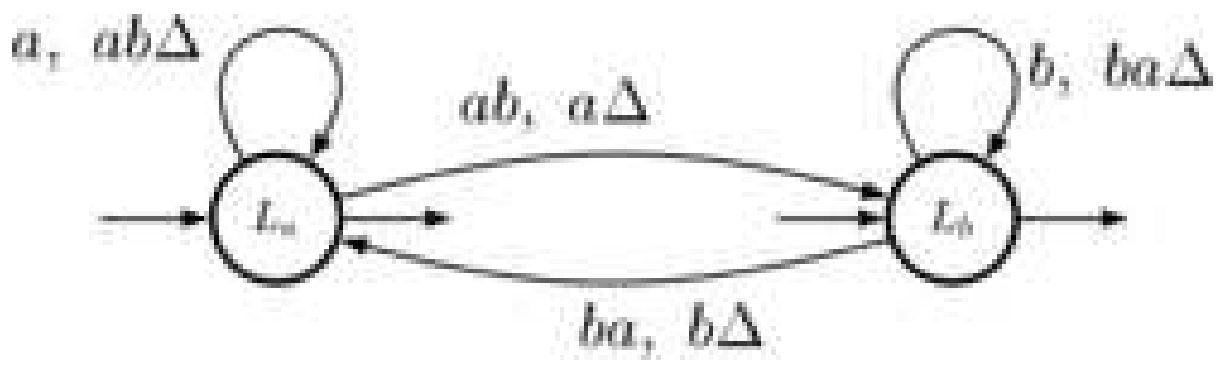

FIG. 7. *Automaton over Σ recognizing the set of all the geodesics $\mathcal{SG}$.*

The automaton of Figure 7 can be viewed as a "symmetrization" of the one of Figure 6. Indeed, the lack of symmetry in the Garside geodesics (due to the handling of $\Delta$) disappears in the set of all the geodesics.

Consider $u = ba * ab * ba * a\Delta \in B_3/Z$. In Figure 5 we have represented with double edges the path corresponding to the Garside normal form of $u$. With simple and plain edges, we have represented the paths corresponding to all the other geodesic representatives of $u$.

Recall that $\pi : \Sigma^* \to B_3/Z$ denotes the canonical homomorphism. Define the map $\psi$ that associates to a geodesic the corresponding Garside geodesic:

$$\psi : \mathcal{SG} \xrightarrow{\pi} B_3/Z \xrightarrow{\phi} \mathcal{G}. \tag{13}$$

The function $\psi$, viewed as a partial function from $\Sigma^*$ to $\Sigma^*$ with domain $\mathcal{SG}$ and image $\mathcal{G}$, is recognized by the transducer $\mathcal{T}$ given in Figure 8. The interpretation of the labels of the states is as follows. Consider a successful path of input label $u \in \mathcal{SG}$. If the path terminates in a state $L_x$, then the last letter of $u$ belongs to $L_x$ and the number of $\Delta$'s in $u$ is even; if the path terminates in a state $(L_x, \Delta)$, then the last letter of $u$ belongs to $L_x$ and the number of $\Delta$'s in $u$ is odd. With this interpretation in mind, proving that $\mathcal{T}$ recognizes $\psi$ is easy.

The transduction $\psi^{-1}$ is recognized by the converse transducer $\mathcal{T}^{-1}$ obtained by exchanging the input and output labels in $\mathcal{T}$. That is,

$$q_1 \xrightarrow{u|v} q_2 \text{ in } \mathcal{T} \iff q_1 \xrightarrow{v|u} q_2 \text{ in } \mathcal{T}^{-1}.$$

The map $\psi^{-1}$ can also be defined recursively as follows. Recall that $\iota$ was defined in (10). We have $\psi^{-1}(1) = 1$, $\forall x \in \Sigma, \psi^{-1}(x) = x$, $\psi^{-1}(\Delta) = \{u \cdot v \in \Sigma^2 | u * v = \Delta\}$, and $\forall u = u_1 \cdots u_k \in \mathcal{G} \backslash \{1, \Delta\}$,

$$\psi^{-1}(u_1 \cdots u_k) = u_1 \cdot \psi^{-1}(u_2 \cdots u_k) \cup u_1\Delta \cdot \psi^{-1}(\iota(u_2) \cdots \iota(u_{k-1})(\iota(u_k)\Delta)).$$

For instance, we have

$$\begin{aligned}\psi^{-1}(a \cdot ab \cdot b) &= a \cdot \psi^{-1}(ab \cdot b) \cup a\Delta \cdot \psi^{-1}(ba \cdot a\Delta) \\ &= a \cdot ab \cdot \psi^{-1}(b) \cup a \cdot ab\Delta \cdot \psi^{-1}(a\Delta) \\ &\quad \cup a\Delta \cdot ba \cdot \psi^{-1}(a\Delta) \cup a\Delta \cdot ba\Delta \cdot \psi^{-1}(b) \\ &= \{a \cdot ab \cdot b\} \cup \{a \cdot ab\Delta \cdot a\Delta\} \\ &\quad \cup \{a\Delta \cdot ba \cdot a\Delta\} \cup \{a\Delta \cdot ba\Delta \cdot b\}.\end{aligned}$$

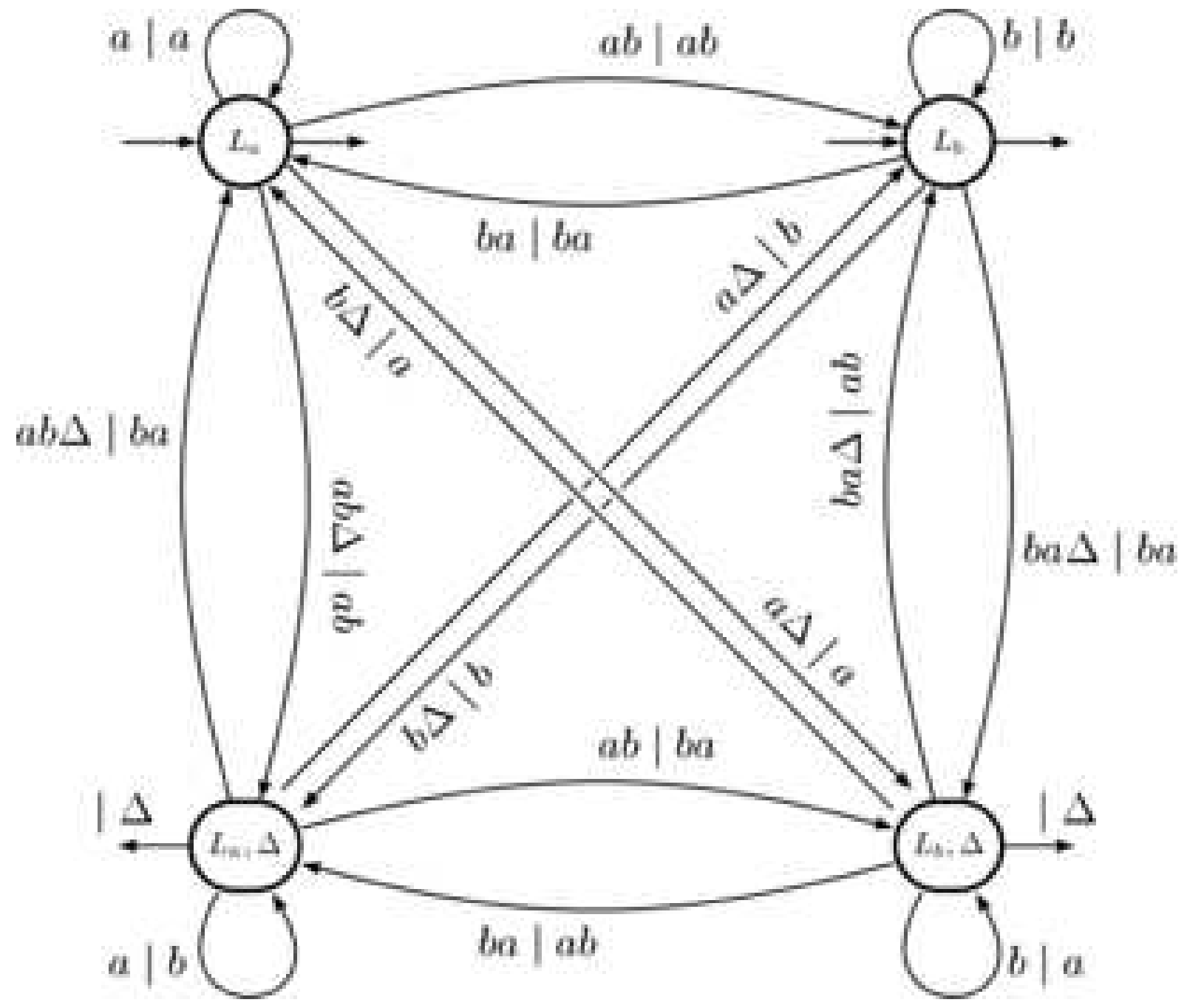

FIG. 8. *The transducer $\mathcal{T}$ recognizing $\psi: \mathcal{SG} \to \mathcal{G}$.*

*Computing lengths in $B_3$.* Proposition 3.2 explains how to go from a Garside normal form in $B_3$ to a geodesic representative with respect to $S = \{a, a^{-1}, b, b^{-1}\}$. See [3, 22] for a proof.

PROPOSITION 3.2. *Let $w = g_1 \cdot \cdots \cdot g_m \cdot \Delta^\delta$ be the Garside normal form of an element $g$ in $B_3$. If $\delta \geq 0$, then $w$ (viewed as a word of $\{a,b\}^*$ by setting $\Delta = aba$) is a minimal length representative of $g$ with respect to $S$. If $\delta < 0$, then a minimal length representative of $g$ with respect to $S$ is obtained as follows. Pick $\min(m, |\delta|)$ elements of $(g_1, \dots, g_m)$ of maximal length with respect to $\{a,b\}$. (This choice may not be unique.) Spread the $\Delta^{-1}$'s in $g_1 \cdot \cdots \cdot g_m$ and, for the picked elements, perform the transformations: $[u\Delta^{-1} \to \widetilde{u}]$, where $\widetilde{a} = b^{-1}a^{-1}, \widetilde{b} = a^{-1}b^{-1}, \widetilde{ab} = b^{-1}, \widetilde{ba} = a^{-1}$.*

The above is best understood on a couple of examples. Below, the left-hand side is a Garside normal form, and the right-hand side is a minimal representative with respect to $S$ of the same element (recall that $a\Delta^{-1} = \Delta^{-1}b, b\Delta^{-1} = \Delta^{-1}a$):

$$a \cdot ab \cdot b \cdot \Delta^{-1} \longrightarrow a \cdot b^{-1} \cdot a,$$

$$a \cdot ab \cdot b \cdot \Delta^{-2} \longrightarrow a \cdot b^{-1} \cdot b^{-1}a^{-1} \text{ or } b^{-1}a^{-1} \cdot a^{-1} \cdot b,$$

$$a \cdot ab \cdot b \cdot \Delta^{-4} \longrightarrow b^{-1}a^{-1} \cdot a^{-1} \cdot a^{-1}b^{-1} \cdot a^{-1}b^{-1}a^{-1}.$$

**4. Random walk on $B_3$.** Consider the random walk $(B_3, \nu)$, where $\nu$ is a probability measure on $p^{-1}(\Sigma) = \{a\Delta^k, ab\Delta^k, b\Delta^k, ba\Delta^k, k \in \mathbb{Z}\}$ such that $\bigcup_n \mathrm{supp}(\nu^{*n}) = B_3$ and $\sum_{x \in p^{-1}(\Sigma)} |x|_S \nu(x) < \infty$. Observe that $S \subset p^{-1}(\Sigma)$, hence, we are in a more general context than the one described in the Introduction.

Let $(x_n)_{n\in\mathbb{N}}$ be a sequence of i.i.d. r.v.'s distributed according to $\nu$. Let $(X_n)_n$ be the corresponding realization of $(B_3, \nu)$. Define $\mu = \nu \circ p^{-1}$ which is a probability distribution on $\Sigma$. The sequence $(p(X_n))_n$ is a realization of the random walk $(B_3/Z, \mu)$. Let $\widehat{X}_n \Delta^{k_n}$ be the Garside normal form of $X_n$. Recall that we have set $S_+ = \{a, b\}$. Define the a.s. limits:

$$\gamma = \lim_{n\to\infty} \frac{|X_n|_S}{n}, \qquad \gamma_\Sigma = \lim_{n\to\infty} \frac{|p(X_n)|_\Sigma}{n},$$
(14)
$$\gamma_{S_+} = \lim_{n\to\infty} \frac{|p(X_n)|_{S_+}}{n}, \qquad \gamma_\Delta = \lim_{n\to\infty} \frac{k_n}{n}.$$

The constant $\gamma$ is the drift on $B_3$, $\gamma_\Sigma$ the drift on $B_3/Z$ with respect to $\Sigma$, $\gamma_{S_+}$ the drift on $B_3/Z$ with respect to $S^+$, and $\gamma_\Delta$ the drift along the center $Z$. Observe that we also have $\gamma_\Sigma = \lim_n |p(X_n)|_T/n = \lim_n |\widehat{X}_n|_T/n$ and $\gamma_{S_+} = \lim_n |\widehat{X}_n|_{S_+}/n$. To obtain $\gamma$, we need to compute $\gamma_\Sigma, \gamma_{S_+}$ and $\gamma_\Delta$. To retrieve $\gamma_\Sigma, \gamma_{S_+}$ and $\gamma_\Delta$, we need to compute the harmonic measure of $(B_3/Z, \mu)$. Besides, determining the harmonic measure is interesting in itself.

It is convenient to view the random walk $(B_3/Z, \mu)$ on the set $\mathcal{G}$ of Garside normal forms. Set $y_n = p(x_n)$ and $Y_n = \phi \circ p(X_n)$. By definition, $(Y_n)_n$ is a Markov chain on $\mathcal{G}$ and satisfies

$$Y_{n+1} = Y_n \circledast y_n = y_0 \circledast y_1 \circledast \cdots \circledast y_n.$$

Specializing (11), we get the following: if $Y_n = 1$, then $Y_{n+1} = y_n$, if $Y_n = \Delta$, then $Y_{n+1} = \iota(y_n)\Delta$, otherwise, if $Y_n = u_1 \cdots u_k \neq \{1, \Delta\}$ and $y_n = x$, then

$$Y_{n+1} = \begin{cases} u_1 \cdots u_{k-1}, & \text{if } x = u_k^{-1}, \\ u_1 \cdots [u_{k-1}\Delta], & \text{if } u_k * x = \Delta, \\ u_1 \cdots u_{k-1}v, & \text{if } u_k * x = v \in \Sigma, \\ u_1 \cdots u_{k-1}u_k x, & \text{if } u_k * x \notin \Sigma \cup \{1, \Delta\} \text{ and } u_k \in T, \\ u_1 \cdots u_{k-1} w_k [\iota(x)\Delta], & \text{if } u_k * x \notin \Sigma \cup \{1, \Delta\} \\ & \text{and } u_k = w_k \Delta \in T^{-1}. \end{cases}$$

In all cases, $Y_n$ and $Y_{n+1}$ differ only in the last two letters. Define the set of infinite words:

$$\mathcal{G}^\infty = \{u_0 u_1 u_2 \cdots \in T^{\mathbb{N}} | \forall i \in \mathbb{N}^* \, \mathrm{Last}(u_i) = \mathrm{First}(u_{i+1})\}. \tag{15}$$

A word belongs to $\mathcal{G}^\infty$ iff all its finite prefixes belong to $\mathcal{G}$. Therefore, $\mathcal{G}^\infty$ can be viewed as the set of infinite Garside geodesics. The set $\mathcal{G}^\infty$ is the set of *ends*, or *end boundary*, of $\mathcal{X}(B_3/Z,\Sigma)$. The group $B_3/Z$ is hyperbolic in the sense of Gromov [11], and $\mathcal{G}^\infty$ coincides with the *hyperbolic boundary* of $B_3/Z$. In the symbolic dynamic terminology, $\mathcal{G}^\infty$ is a one-sided subshift of finite type, since it is defined through a finite set of forbidden patterns.

Equip $T^{\mathbb{N}}$ with the Borel $\sigma$-algebra associated with the product topology. Denote by $(u_1\cdots u_n T^{\mathbb{N}})$ the order-$n$ cylinder in $T^{\mathbb{N}}$ defined by $u_1\cdots u_n$.

Recall that $(Y_n)_n$ is transient and that $Y_n$ and $Y_{n+1}$ differ only in the last two letters. It follows easily that there exists a r.v. $Y_\infty$ such that a.s.

$$\lim_{n\to\infty} Y_n = Y_\infty,$$

meaning that the length of the common prefix between $Y_n$ and $Y_\infty$ goes to infinity. Let $\mu^\infty$ be the distribution of $Y_\infty$. The measure $\mu^\infty$ captures the direction in which the random walk $(B_3/Z,\mu)$ goes to infinity. We call $\mu^\infty$ the *harmonic measure* of $(B_3/Z,\mu)$. The pair $(\mathcal{G}^\infty,\mu^\infty)$ is the *Poisson boundary* of $(B_3/Z,\mu)$, and $\mathcal{G}^\infty$ is the *Martin boundary*, as well as the *minimal Martin boundary* of $(B_3/Z,\mu)$. See [13, 14] for details.

The harmonic measure contains a lot of information about the random walk. For instance, one can easily retrieve the drift $\gamma_\Sigma$ knowing $\mu^\infty$; see Proposition 4.1.

The map $\mathcal{G}\times\mathcal{G}\to\mathcal{G}, (u,v)\mapsto u\circledast v$, defined in (9), extends naturally to a group action $\mathcal{G}\times\mathcal{G}^\infty\to\mathcal{G}^\infty, (u,\xi)\mapsto u\circledast\xi$. Below we focus on its restriction:

$$\Sigma\times\mathcal{G}^\infty\to\mathcal{G}^\infty, \qquad (u,\xi)\mapsto u\circledast\xi. \tag{16}$$

Given a probability measure $p^\infty$ on $\mathcal{G}^\infty$, and given $u\in\mathcal{G}$, define the measure $u\circledast p^\infty$ by $\forall f:\mathcal{G}^\infty\to\mathbb{R}, \int f(\xi)\,d(u\circledast p^\infty)(\xi)=\int f(u\circledast\xi)\,dp^\infty(\xi)$. A probability measure $p^\infty$ on $\mathcal{G}^\infty$ is *$\mu$-stationary* if

$$p^\infty(\cdot)=\sum_{u\in\Sigma}\mu(u)[u\circledast p^\infty](\cdot). \tag{17}$$

PROPOSITION 4.1. *The probability $\mu^\infty$ is the unique $\mu$-stationary probability measure on $\mathcal{G}^\infty$. For $u\in\Sigma$, let $f(u)$ be the unique element $v$ of $T$ such that $u*v\in\{1,\Delta\}$. We have a.s.*

$$\gamma_\Sigma=\sum_{u\in\Sigma}\mu(u)\left[-\mu^\infty(f(u)T^{\mathbb{N}})+\sum_{\substack{v\in T\\ u*v\notin\Sigma\cup\{1,\Delta\}}}\mu^\infty(vT^{\mathbb{N}})\right]. \tag{18}$$

The result on $\mu^\infty$ is classical. It is proved, for instance, in [14], Section 7, in the more general context of random walks on Gromov hyperbolic groups.

4.1. *Statements of the results.* Given $u \in \Sigma$, set $\mathrm{Next}(u) = \{v \in \Sigma | \, \mathrm{Last}(u) = \mathrm{First}(v)\}$. We have

$$\mathrm{Next}(a) = \mathrm{Next}(ba) = \mathrm{Next}(b\Delta) = \mathrm{Next}(ab\Delta) = \{a, ab, a\Delta, ab\Delta\},$$
$$\mathrm{Next}(b) = \mathrm{Next}(ab) = \mathrm{Next}(a\Delta) = \mathrm{Next}(ba\Delta) = \{b, ba, b\Delta, ba\Delta\}.$$

Observe that we have $\mathcal{G}^\infty = \{u_0 u_1 \cdots \in T^{\mathbb{N}} | \forall i \ u_{i+1} \in \mathrm{Next}(u_i)\}$. The *traffic equations* of $(B_3/Z, \mu)$ are the equations in the indeterminates $x(u), u \in \Sigma$, defined by, $\forall u \in \Sigma$,

$$\begin{aligned} x(u) = \mu(u) \sum_{v \in \mathrm{Next}(u)} x(v) + \sum_{\substack{v,w \in \Sigma \\ v * w = u}} \mu(v) x(w) \\ + \sum_{\substack{v,w \in \Sigma, v*w=1 \\ u \in \mathrm{Next}(w)}} \mu(v) \frac{x(w)}{\sum_{z \in \mathrm{Next}(w)} x(z)} x(u) \\ + \sum_{\substack{v,w \in \Sigma, v*w=\Delta \\ \iota(u)\Delta \in \mathrm{Next}(w)}} \mu(v) \frac{x(w)}{\sum_{z \in \mathrm{Next}(w)} x(z)} x(\iota(u)\Delta). \end{aligned} \tag{19}$$

THEOREM 4.2. *Consider the random walk $(B_3, \nu)$, where $\nu$ is a probability measure on $p^{-1}(\Sigma) = \{a\Delta^k, ab\Delta^k, b\Delta^k, ba\Delta^k, k \in \mathbb{Z}\}$ such that $\bigcup_n \mathrm{supp}(\nu^{*n}) = B_3$ and $\sum_x |x|_S \nu(x) < \infty$. Set $\mu = \nu \circ p^{-1}$. The traffic equations of $(B_3/Z, \mu)$ have a unique solution $r \in \{x \in (\mathbb{R}_+^*)^\Sigma | \sum_{u \in \Sigma} x(u) = 1\}$. For $u \in \Sigma$, set*

$$q(u) = r(u)/r(\mathrm{Next}(u)), \qquad R(u) = r(u) + r(u\Delta).$$

*The harmonic measure $\mu^\infty$ is given by, $\forall u = u_1 \cdots u_k \in \mathcal{G} \cap T^k$,*

$$\mu^\infty(uT^{\mathbb{N}}) = \sum_{v_1 \cdots v_k \in \psi^{-1}(u)} q(v_1) \cdots q(v_{k-1}) R(v_k). \tag{20}$$

Here is an equivalent way of stating (20). Consider the automaton over the alphabet $T$ with multiplicities over $(\mathbb{R}_+, +, \times)$ represented in Figure 9.

The automaton of Figure 9 has the same structure as the one of Figure 8, and a label $[u|v]$ in Figure 8 has been replaced by a label $[v|q(u)]$ in Figure 9. Let $Q = (L_a; [L_b, \Delta]; L_b; [La, \Delta])$ be the states. Let $\alpha \in \mathbb{R}_+^{1 \times Q}, \alpha = [1,0,1,0]$, be the initial vector, and let $\mathcal{M} : T \to \mathbb{R}_+^{Q \times Q}$ be the matrices associated with the automaton. Define $\beta : T \to \mathbb{R}_+^{Q \times 1}$ by, if $u = a$ or $ab$, $\beta(u) = [R(u), R(\iota(u)), 0, 0]^T$; if $u = b$ or $ba$, $\beta(u) = [0, 0, R(u), R(\iota(u))]^T$. The measure $\mu^\infty$ is the *rational* measure defined by

$$\forall k, \forall u = u_1 \cdots u_k \in T^k \qquad \mu^\infty(uT^{\mathbb{N}}) = \alpha \mathcal{M}(u_1) \cdots \mathcal{M}(u_{k-1}) \beta(u_k). \tag{21}$$

In particular, if $u \notin \mathcal{G}$, we check that $\alpha \mathcal{M}(u_1) \cdots \mathcal{M}(u_{k-1}) \beta(u_k) = 0$.

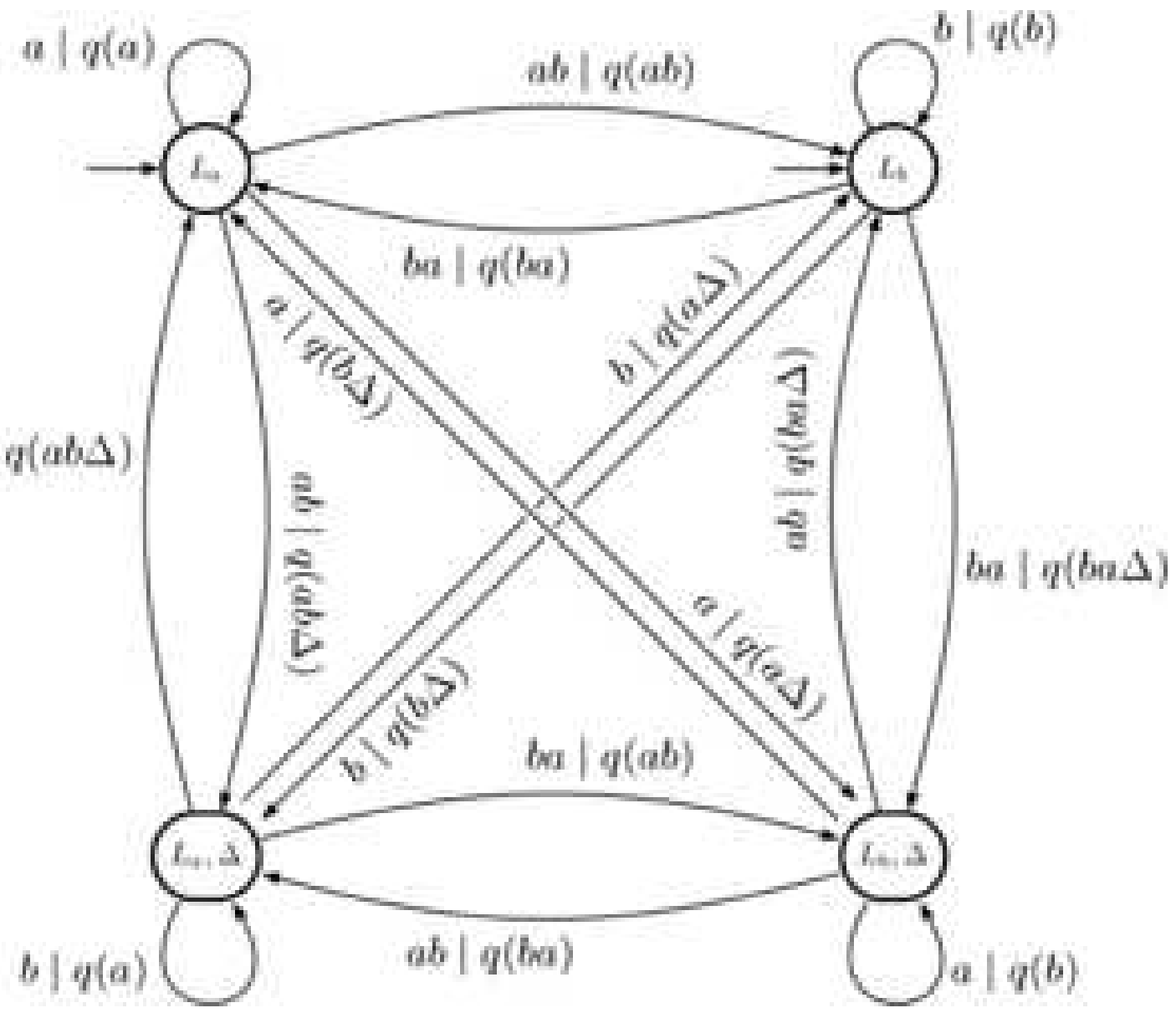

FIG. 9. *The $\mathbb{R}_+$-automaton defining the harmonic measure $\mu^\infty$.*

Let $\ell : \Sigma \to T$ be defined by $\ell(u) = u$ if $u \in T$, and $\ell(u) = v$ if $u = v\Delta \in T^{-1}$. Define $P \in \mathbb{R}_+^{\Sigma\times\Sigma}$ by

$$P_{u,v} = \begin{cases} r(v)/r(\mathrm{Next}(u)), & \text{if } v \in \mathrm{Next}(u), \\ 0, & \text{otherwise.} \end{cases} \tag{22}$$

By the Perron–Frobenius theorem, the equations $xP = x$ have a unique solution $\mathfrak{p} \in \{x \in (\mathbb{R}_+^*)^\Sigma \mid \sum_{u\in\Sigma} x(u) = 1\}$. Set $\mathfrak{p}(1) = \mathfrak{p}(a) + \mathfrak{p}(a\Delta) + \mathfrak{p}(b) + \mathfrak{p}(b\Delta)$ and $\mathfrak{p}(2) = 1 - \mathfrak{p}(1)$. We have the following:

THEOREM 4.3. *The drifts $\gamma_\Sigma$, $\gamma_{S_+}$ and $\gamma_\Delta$, defined in (14), satisfy*

$$\gamma_\Sigma = \sum_{u\in T}\left[-\mu(u^{-1}) - \mu(\Delta u^{-1}) + \sum_{v\mid u\in\mathrm{Next}(v)} \mu(v)\right] R(u), \tag{23}$$

$$\begin{aligned} \gamma_{S_+} = \sum_{u\in T}\Bigg[&-\mu(u^{-1})|u|_{S^+} - \mu(\Delta u^{-1})|u|_{S^+} \\ &+ \sum_{v*u\in\Sigma} \mu(v)(|\ell(v*u)|_{S^+} - |u|_{S^+}) \\ &+ \sum_{v\mid u\in\mathrm{Next}(v)} \mu(v)|v|_{S^+}\Bigg] R(u), \end{aligned} \tag{24}$$

$$(25)\ \gamma_\Delta = \sum_{u\in T}\left[\sum_{v*u\in T^{-1}\cup\Delta}\ \sum_{k\in\mathbb{Z}}(k+1)\nu(\Delta^k v) + \sum_{v*u\notin T^{-1}\cup\Delta}\ \sum_{k\in\mathbb{Z}} k\nu(\Delta^k v)\right]R(u).$$

*The following relation holds:*

$$\gamma_{S_+} = [\mathfrak{p}(1)+2\mathfrak{p}(2)]\gamma_\Sigma. \tag{26}$$

*The drift of the random walk* $(B_3,\nu)$ *is then given by*

$$\gamma = \begin{cases} \gamma_{S_+} + 3\gamma_\Delta, & \text{if } \gamma_\Delta \geq 0,\\ \gamma_{S_+} - |\gamma_\Delta|, & \text{if } \gamma_\Delta < 0, |\gamma_\Delta| \leq \mathfrak{p}(2)\gamma_\Sigma,\\ \gamma_{S_+} + |\gamma_\Delta| - 2\mathfrak{p}(2)\gamma_\Sigma, & \text{if } \gamma_\Delta < 0, \mathfrak{p}(2)\gamma_\Sigma < |\gamma_\Delta| \leq \gamma_\Sigma,\\ -\gamma_{S_+} + 3|\gamma_\Delta|, & \text{if } \gamma_\Delta < 0, \gamma_\Sigma < |\gamma_\Delta|. \end{cases} \tag{27}$$

Now let us restrict ourselves to the case considered in Figure 2: the measure $\nu$ is supported by $S$ and $p = \nu(a) = 1/2 - \nu(b^{-1}), q = \nu(b) = 1/2 - \nu(a^{-1})$. Then the drifts take a simpler form. Recall that, for $u \in T$, $R(u) = r(u) + r(u\Delta)$. We have

$$\begin{aligned} \gamma_\Sigma &= qR(a) + pR(b) + (1/2-q)R(ab) + (1/2-p)R(ba),\\ \gamma_{S_+} &= 1/2 - p - q + 3qR(a) + 3pR(b),\\ \gamma_\Delta &= p + q - 1/2 - qR(a) - pR(b). \end{aligned} \tag{28}$$

If $p+q \geq 1/2$, we are in one of the first two cases in (27) for $\gamma(p,q)$. And, more precisely,

$$\gamma = \begin{cases} 2(p+q)-1, & \text{if } \gamma_\Delta \geq 0,\\ 2qR(a)+2pR(b), & \text{if } \gamma_\Delta < 0. \end{cases} \tag{29}$$

If $p+q < 1/2$, then $\gamma(p,q) = \gamma(1/2-p, 1/2-q)$.

The three lines of nondifferentiability in Figure 2 can be identified on the above equations. One corresponds to $p+q = 1/2$. Another one corresponds to the separation between the two cases $[p+q \geq 1/2, \gamma_\Delta \geq 0]$ and $[p+q \geq 1/2, \gamma_\Delta < 0]$ in (29); and the last one is the reflection of the second one on $p+q < 1/2$.

Here are the traffic equations in the case $p = \nu(a) = 1/2 - \nu(b^{-1}), q = \nu(b) = 1/2 - \nu(a^{-1})$. Set $S_a = r(a) + r(ab) + r(a\Delta) + r(ab\Delta)$ and $S_b = r(b) + r(ba) + r(b\Delta) + r(ba\Delta)$. We have

$$\begin{aligned} r(a) &= pS_a + q(r(ab\Delta)r(a)S_a^{-1} + r(ab)r(b\Delta)S_b^{-1})\\ &\qquad + (1/2-q)(r(a)^2S_a^{-1} + r(a\Delta)r(b\Delta)S_b^{-1}) + (1/2-p)r(ba),\\ r(b) &= qS_b + p(r(ba\Delta)r(b)S_b^{-1} + r(ba)r(a\Delta)S_a^{-1})\\ &\qquad + (1/2-p)(r(b)^2S_b^{-1} + r(b\Delta)r(a\Delta)S_a^{-1}) + (1/2-q)r(ab),\\ r(ab) &= pr(b) + q(r(ab\Delta)r(ab)S_a^{-1} + r(ab)r(ba\Delta)S_b^{-1}) \end{aligned}$$

$$
\begin{aligned}
&\quad + (1/2-q)(r(a)r(ab)S_a^{-1} + r(a\Delta)r(ba\Delta)S_b^{-1}),\\
r(ba) &= qr(a) + p(r(ba\Delta)r(ba)S_b^{-1} + r(ba)r(ab\Delta)S_a^{-1})\\
&\quad + (1/2-p)(r(b)r(ba)S_b^{-1} + r(b\Delta)r(ab\Delta)S_a^{-1}),\\
r(a\Delta) &= q(r(ab\Delta)r(a\Delta)S_a^{-1} + r(ab)r(b)S_b^{-1})\\
&\quad + (1/2-q)(r(a)r(a\Delta)S_a^{-1} + r(a\Delta)r(b)S_b^{-1}) + (1/2-p)r(ba\Delta),\\
r(b\Delta) &= p(r(ba\Delta)r(b\Delta)S_b^{-1} + r(ba)r(a)S_a^{-1})\\
&\quad + (1/2-p)(r(b)r(b\Delta)S_b^{-1} + r(b\Delta)r(a)S_a^{-1}) + (1/2-q)r(ab\Delta),\\
r(ab\Delta) &= pr(b\Delta) + q(r(ab)r(ba)S_b^{-1} + r(ab\Delta)^2S_a^{-1})\\
&\quad + (1/2-p)S_a + (1/2-q)(r(a)r(ab\Delta)S_a^{-1} + r(a\Delta)r(ba)S_b^{-1}),\\
r(ba\Delta) &= qr(a\Delta) + p(r(ba)r(ab)S_a^{-1} + r(ba\Delta)^2S_b^{-1})\\
&\quad + (1/2-q)S_b + (1/2-p)(r(b)r(ba\Delta)S_b^{-1} + r(b\Delta)r(ab)S_a^{-1}).
\end{aligned}
\tag{30}
$$

Under additional symmetries on $\nu$, closed form formulas for the drifts are obtained.

Proposition 4.4. *Assume that $\nu(a) = \nu(a^{-1}) = p$, $\nu(b) = \nu(b^{-1}) = 1/2 - p$, with $p \in (0, 1/4]$. Let $u$ be the smallest root in $(0, 1)$ of the polynomial*

$$P = 2(4p-1)X^3 + (24p^2 - 18p + 1)X^2 + p(-12p + 7)X + p(2p-1).$$

*The various drifts are given by*

$$\gamma(p) = \gamma_\Sigma(p) = -2\gamma_\Delta(p) = \tfrac{2}{3}\gamma_{S_+}(p) = p + (1-4p)u.$$

*For $p \in [1/4, 1/2)$, we have $\gamma(p) = \gamma(1/2-p)$ and similarly for $\gamma_\Sigma, \gamma_{S_+}$ and $\gamma_\Delta$.*

Proposition 4.5. *Assume that $\nu(a) = \nu(b) = p$, $\nu(a^{-1}) = \nu(b^{-1}) = 1/2 - p$, with $p \in (0, 1/2)$. We have*

$$
\begin{aligned}
\gamma_\Sigma(p) &= \frac{-1 + \sqrt{16p^2 - 8p + 5}}{4},\\
\gamma_{S^+}(p) &= \frac{4p^2 + p + 1 - 3p\sqrt{16p^2 - 8p + 5}}{2(1-4p)},\\
\gamma_\Delta(p) &= \frac{12p^2 - 5p + 1 - p\sqrt{16p^2 - 8p + 5}}{2(1-4p)}.
\end{aligned}
$$

*And, eventually,*

$$\gamma(p) = \max\Big[1-4p, \frac{(1-2p)(-1-4p+\sqrt{5-8p+16p^2})}{2(1-4p)}, \frac{p(-3+4p+\sqrt{5-8p+16p^2})}{-1+4p}, -1+4p\Big].$$

In both cases, the method is to simplify the traffic equations (30) by using the additional symmetries. The resulting equations can be completely solved, leading to the above results. The details of the computations are given in [23], Sections 2 and 3.

REMARK 4.6. One can retrieve from the harmonic measure $\mu^\infty$ other quantities of interest for the random walk $(B_3/Z, \mu)$: (a) the entropy; (b) the minimal positive harmonic functions; or (c) the Green function. Details are provided in [23], Section 5.

REMARK 4.7. Consider the sets $\overline{\Sigma} = \Sigma \cup \{1, \Delta\}$ and $p^{-1}(\overline{\Sigma}) = \{\Delta^k, a\Delta^k, ab\Delta^k, b\Delta^k, ba\Delta^k, k \in \mathbb{Z}\}$. The above results extend to the random walks $(B_3, \nu)$ and $(B_3/Z, \mu)$, where $\nu$ and $\mu$ are probabilities generating the whole group and with respective support $p^{-1}(\overline{\Sigma})$ and $\overline{\Sigma}$. The traffic equations of $(B_3/Z, \mu)$ are still over the indeterminates $(x(u))_{u\in\Sigma}$. They are given by (19) with $[v \in \Sigma]$ being replaced by $[v \in \overline{\Sigma}]$.

4.2. *Proof of Theorem 4.2.* Let us provide a sketch of the proof. First we introduce an auxiliary process $(Z_n)_n$ evolving on the set $\mathcal{SG}$, and such that $\psi(Z_n) = Y_n$. We show that the limit $Z_\infty = \lim_n Z_n$ exists, and that its law has a simple Markovian structure. Then the formulas (20) for $\mu^\infty$ and (23) for $\gamma_\Sigma$ are obtained using Proposition 4.1. Second, the drift along the center $\gamma_\Delta$ is computed. Third, we compute $\gamma$ from $\gamma_\Sigma$ and $\gamma_\Delta$, using Proposition 3.2.

Define the map

$$\begin{aligned} \Sigma \times \mathcal{SG} &\to \mathcal{SG}, \\ (u, v) &\mapsto u \odot v \end{aligned} \tag{31}$$

as follows: $u \odot 1 = u$, $u \odot \Delta = u\Delta$, and for $v = v_1 \cdots v_k \in \mathcal{SG}\backslash\{1, \Delta\}$,

$$u \odot v = \begin{cases} v_2 \cdot v_3 \cdots v_k, & \text{if } u = v_1^{-1}, \\ [\iota(v_2)\Delta] \cdot v_3 \cdots v_k, & \text{if } u * v_1 = \Delta, \\ w \cdot v_2 \cdot v_3 \cdots v_k, & \text{if } u * v_1 = w \in \Sigma, \\ u \cdot v_1 \cdot v_2 \cdots v_k, & \text{if } u * v_1 \notin \Sigma \cup \{1, \Delta\}. \end{cases} \tag{32}$$

By construction, if $v$ satisfies the criterion in (12), then $u \odot v$ also satisfies it. This shows that we have indeed a map $\Sigma \times \mathcal{SG} \to \mathcal{SG}$. Given $x, y \in \Sigma, x * y = z \in \Sigma$, we do not have, in general, $x \odot [y \odot v] = z \odot v$. For instance,

$$a \odot [b \odot (ab \cdot b)] = a \cdot a\Delta, \qquad ab \odot (ab \cdot b) = a\Delta \cdot b.$$

Extend the sequence $(y_n)_{n \in \mathbb{N}}$ to a bi-infinite i.i.d. sequence $(y_n)_{n \in \mathbb{Z}}$. Let $(\Omega, \mathcal{F}, P)$ be the underlying probability space. Consider that we have, $\forall n \in \mathbb{Z}$, $y_n = y_0 \circ \theta^n$, where $\theta : \Omega \to \Omega$ is a bijective and measure-preserving shift.

Define the random process $(Z_n)_{n \in \mathbb{N}}$ on the state space $\mathcal{SG}$ by $Z_0 = 1$ and

$$\begin{aligned} Z_{n+1} &= y_0 \odot [Z_n \circ \theta] \\ &= y_0 \odot [y_1 \odot [y_2 \odot [\cdots [y_{n-1} \odot y_n]] \cdots ]]. \end{aligned} \tag{33}$$

By definition, we have $\pi(Y_n) = y_0 * \cdots * y_{n-1}$ and $\pi(Z_n) = y_0 * \cdots * y_{n-1}$. Since $Z_n$ belongs to $\mathcal{SG}$ and $Y_n$ belongs to $\mathcal{G}$, we conclude that $\psi(Z_n) = Y_n$.

The process $(W_n)_{n \in \mathbb{N}} = (Z_n \circ \theta^{-n})_{n \in \mathbb{N}}$ is Markovian and satisfies $W_{n+1} = y_{-(n+1)} \odot W_n$. For every $n$, the law of $Z_n$ is the same as the law of $W_n$, but the law of the process $(Z_n)_n$ is different from the law of the process $(W_n)_n$. In particular, $(Z_n)_n$ is not Markovian. Considering the process $(Z_n)_n$ amounts to performing a *coupling from the past* construction [24], also known as a Loynes backward construction in queueing theory [16].

Let us compare the trajectories of $(Y_n)_n$ and $(Z_n)_n$ on an example. Consider the driving sequence $y_0 = a, y_1 = ab\Delta, y_2 = a, y_3 = ab\Delta$. We have

$$\begin{aligned} &Y_1 = a, \qquad Y_2 = a \cdot ab\Delta, \qquad Y_3 = a \cdot ab \cdot b\Delta, \qquad Y_4 = a \cdot ab \cdot b \cdot ba, \\ &Z_1 = a, \qquad Z_2 = a \cdot ab\Delta, \qquad Z_3 = a \cdot ab\Delta \cdot a, \qquad Z_4 = a \cdot ab\Delta \cdot a \cdot ab\Delta. \end{aligned}$$

For $u, v \in \Sigma^*$, let $u \wedge v \in \Sigma^*$ be the longest common prefix of $u$ and $v$. Set $d(u, v) = \max(|u|, |v|) - |u \wedge v|$. This defines a distance on $\Sigma^*$ called the *prefix distance*.

LEMMA 4.8. *We have $d(Z_n, Z_{n+1}) \leq 2$, where $d(\cdot, \cdot)$ is the prefix distance on $\Sigma^*$.*

PROOF. If $|Z_n|_\Sigma = 0$ or 1, then the result is immediate. Assume that $|Z_n|_\Sigma \geq 2$. Let $u, v$ be the $\Sigma$-valued r.v.'s such that $Z_n \in \Sigma^* uv$. For $u_1, \ldots, u_n \in \Sigma$, set

$$\psi(u_1, \ldots, u_n) = u_1 \odot [u_2 \odot [u_3 \odot [\cdots [u_{n-1} \odot u_n]] \cdots ]].$$

Define $m = \max\{k \geq 0 | \forall 0 \leq j \leq k \ \ \psi(y_j, y_{j+1}, \ldots, y_{n-1}) \in \Sigma^* uv\}$. By definition of $\odot$, we must have $\psi(y_m, \ldots, y_{n-1}) = uv$ and

$$\begin{aligned} Z_n &= \psi(y_0, \ldots, y_{m-1}) \psi(y_m, \ldots, y_{n-1}), \\ Z_{n+1} &= \psi(y_0, \ldots, y_{m-1}) \psi(y_m, \ldots, y_n). \end{aligned}$$

We conclude easily. □

Define the set [compare with (15)]

$$\mathcal{SG}^\infty = \{u_0u_1u_2\cdots \in \Sigma^{\mathbb{N}} | \forall i \in \mathbb{N}^* \ \ \mathrm{Last}(u_i) = \mathrm{First}(u_{i+1})\}. \tag{34}$$

The set $\mathcal{SG}^\infty$ is the set of all infinite geodesics.

Since $\psi(Z_n) = Y_n$, the process $(Z_n)_n$ is transient. We deduce from Lemma 4.8 that $\lim_n Z_n = Z_\infty$, where $Z_\infty$ is a $\mathcal{SG}^\infty$-valued r.v. Let us denote by $\kappa^\infty$ the law of $Z_\infty$. We have

$$\forall u \in T^* \cap \mathcal{G} \qquad \mu^\infty(uT^{\mathbb{N}}) = \sum_{v\in\psi^{-1}(u)} \kappa^\infty(v\Sigma^{\mathbb{N}}) + \sum_{v\in\psi^{-1}(u\Delta)} \kappa^\infty(v\Sigma^{\mathbb{N}}). \tag{35}$$

The map $\Sigma \times \mathcal{SG} \to \mathcal{SG}, (u,v) \mapsto u \odot v$, defined in (31), extends naturally to a map

$$\Sigma \times \mathcal{SG}^\infty \to \mathcal{SG}^\infty, \qquad (u,\xi) \mapsto u \odot \xi. \tag{36}$$

The two maps defined in (16) and (36) have very different behaviors: the map in (36) has a local action on the infinite word, as opposed to the map in (16). Consider, for instance, the infinite word $\xi = a \cdot a \cdot a \cdot a \cdots$, which belongs both to $\mathcal{G}^\infty$ and $\mathcal{SG}^\infty$. We have

$$ab \circledast \xi = b \cdot b \cdot b \cdot b \cdots, \qquad ab \odot \xi = b\Delta \cdot a \cdot a \cdot a \cdots.$$

This can be viewed as the reason why it is crucial to work with $(Z_n)_n$ rather than $(Y_n)_n$.

Given a probability measure $p^\infty$ on $\mathcal{SG}^\infty$, and given $u \in \Sigma$, define the measure $u \odot p^\infty$ by $\forall f : \mathcal{SG}^\infty \to \mathbb{R}, \int f(\xi)\, d(u \odot p^\infty)(\xi) = \int f(u \odot \xi)\, dp^\infty(\xi)$. A probability measure $p^\infty$ on $\mathcal{SG}^\infty$ is *$\mu$-stationary* if the analog of (17) (with $\odot$ instead of $\circledast$) holds.

Following step by step (and slightly adapting) the proof of [15], Theorem 1.12, we get the following analog of Proposition 4.1:

LEMMA 4.9. *The probability $\kappa^\infty$ is the unique $\mu$-stationary probability measure on $\mathcal{SG}^\infty$.*

Define $\mathcal{B} = \{x \in (\mathbb{R}_+^*)^\Sigma | \sum_{u\in\Sigma} x(u) = 1\}$. Consider $r \in \mathcal{B}$. The *Markovian multiplicative measure associated with* $r$ is the probability measure $p^\infty$ on $\mathcal{SG}^\infty$ defined by

$$\forall u_1 \cdots u_k \in \mathcal{SG} \qquad p^\infty(u_1 \cdots u_k \Sigma^{\mathbb{N}}) = q(u_1) \cdots q(u_{k-1}) r(u_k), \tag{37}$$

where $q(u) = r(u)/r(\mathrm{Next}(u))$.

The next step consists in showing the following: (i) the traffic equations (19) have a unique solution $r$ in $\mathcal{B}$; (ii) the measure $\kappa^\infty$ is the Markovian multiplicative measure associated with $r$. The proof of (i) and (ii) follows the lines of Lemma 4.4 and Theorem 4.5 in [19]. We now detail the argument.

Consider the probability distribution $\kappa^\infty$ of $Z_\infty$. Assume that $\kappa^\infty$ is the Markovian multiplicative measure associated with some $r \in \mathcal{B}$. According to Lemma 4.9, the measure $\kappa^\infty$ is characterized by the fact that it is $\mu$-stationary. By definition, $\kappa^\infty$ is $\mu$-stationary iff $\forall u = u_1 \cdots u_k \in \mathcal{SG}$,

$$\begin{aligned}\kappa^\infty(u\Sigma^{\mathbb{N}}) = \mu(u_1)\kappa^\infty(u_2\cdots u_k\Sigma^{\mathbb{N}}) &+ \sum_{v*w=u_1} \mu(v)\kappa^\infty(wu_2\cdots u_k\Sigma^{\mathbb{N}}) \\ &+ \sum_{v*w=1} \mu(v)\kappa^\infty(wu_1u_2\cdots u_k\Sigma^{\mathbb{N}}) \\ &+ \sum_{v*w=\Delta} \mu(v)\kappa^\infty(w[\iota(u_1)\Delta]u_2\cdots u_k\Sigma^{\mathbb{N}}).\end{aligned} \tag{38}$$

We can simplify the above using (37), the equation that we get depends on $u$ only via $u_1$:

$$\begin{aligned}r(u_1) = \mu(u_1)\sum_{v\in \mathrm{Next}(u_1)} r(v) &+ \sum_{v*w=u_1}\mu(v)r(w) \\ &+ \sum_{\substack{v*w=1\\ u_1\in\mathrm{Next}(w)}} \mu(v)\frac{r(w)}{\sum_{z\in\mathrm{Next}(w)} r(z)} r(u_1) \\ &+ \sum_{\substack{v*w=\Delta\\ \iota(u_1)\Delta\in\mathrm{Next}(w)}} \mu(v)\frac{r(w)}{\sum_{z\in\mathrm{Next}(w)} r(z)} r(\iota(u_1)\Delta).\end{aligned} \tag{39}$$

To be able to simplify (38), we use the following consequences of (8):

$$\begin{aligned} v*w=u_1 \quad \Longrightarrow \quad & \mathrm{Next}(w) = \mathrm{Next}(u_1) \text{ and also} \\ & \mathrm{Next}(u_1) = \mathrm{Next}(\iota(u_1)\Delta).\end{aligned} \tag{40}$$

The equations (39) are precisely the traffic equations (19) of $(B_3/Z, \mu)$. To summarize, $\kappa^\infty$ is the Markovian multiplicative measure associated with $r \in \mathcal{B}$ iff $r \in \mathcal{B}$ is a solution to the traffic equations (19).

Remark. Assume that $\mu^\infty$ is Markovian multiplicative associated with some $r$. Now write the analog of (38) and try to simplify using (37). Full simplification does not occur anymore and we get a different set of equations for $r$ for each different word $u$. Therefore, the argument fails. Indeed, $\mu^\infty$ is *not* Markovian multiplicative; see (20).

A consequence of the above is that the traffic equations (19) has at most one solution in $\mathcal{B}$. Now, let us prove that they always have a solution in $\mathcal{B}$.

In equation (19), denote the right-hand side of the equality by $\Phi(x)(u)$. This defines an application $\Phi : (\mathbb{R}_+^*)^\Sigma \longrightarrow (\mathbb{R}_+^*)^\Sigma$. The traffic equations are the equations $x = \Phi(x)$.

Consider $X \in \mathcal{B}$. By rearranging the terms in the sum, we check easily that

$$\sum_{u\in\Sigma} \Phi(X)(u) = \sum_{u\in\Sigma} \mu(u) = 1.$$

Hence, we have $\Phi(\mathcal{B}) \subset \mathcal{B}$. Define $\overline{\mathcal{B}} = \{x \in \mathbb{R}_+^\Sigma | \sum_u x(u) = 1\}$, which is the closure of $\mathcal{B}$ and a convex compact subset of $\mathbb{R}^\Sigma$. Observe that $\Phi$ cannot be extended continuously on $\overline{\mathcal{B}}$. For $x \in \overline{\mathcal{B}} \backslash \mathcal{B}$, let $\Phi(x) \subset \overline{\mathcal{B}}$ be the set of possible limits of $\Phi(x_n), x_n \in \mathcal{B}, x_n \to x$. We have extended $\Phi$ to a correspondence $\Phi : \overline{\mathcal{B}} \twoheadrightarrow \overline{\mathcal{B}}$. Clearly, this correspondence has a closed graph and nonempty convex values. Therefore, we are in the domain of application of the Kakutani–Fan–Glicksberg theorem; see [1], Chapter 16. The correspondence has at least one fixed point: $\exists x \in \overline{\mathcal{B}}$ such that $x \in \Phi(x)$. Now, using the shape of the traffic equations, we obtain easily that $x \in \mathcal{B}$.

So we have proved that the traffic equations have exactly one solution in $\mathcal{B}$. The measure $\kappa^\infty$ is the Markovian multiplicative measure associated with this solution.

Since $\mu^\infty$ and $\kappa^\infty$ are linked by (35), this completes the proof of Theorem 4.2.

4.3. *Proof of Theorem* 4.3. First, using Proposition 4.1 and Theorem 4.2, we get the expression for $\gamma_\Sigma$ given in (23).

Now consider $\gamma_{S+}$. Recall that $\ell : \Sigma \to T$ is defined by $\ell(u) = u$ if $u \in T$, and $\ell(u) = v$ if $u = v\Delta \in T^{-1}$. Using the same argumentation as for (18) in Proposition 4.1, we have

$$\begin{aligned}\gamma_{S_+} = \sum_{u\in T}\Bigg[&-\mu(u^{-1})|u|_{S^+} - \mu(\Delta u^{-1})|u|_{S^+}\\ &+ \sum_{v*u\in\Sigma} \mu(v)(|\ell(v*u)|_{S^+} - |u|_{S^+})\\ &+ \sum_{v|u\in\mathrm{Next}(v)} \mu(v)|v|_{S^+}\Bigg]\mu^\infty(uT^{\mathbb{N}}).\end{aligned}$$

Since we just proved that $\mu^\infty(uT^{\mathbb{N}}) = R(u)$, we get (24).

Let us focus on $\gamma_\Delta$. Set $U = p^{-1}(\Sigma) = \{a\Delta^k, ab\Delta^k, b\Delta^k, ba\Delta^k, k \in \mathbb{Z}\}$. Define the function $\theta_\Delta : U^{\mathbb{N}} \times \mathcal{G}^\infty \to \mathbb{Z}$ in the following way: for all $(\omega, \xi) \in$

$U^{\mathbb{N}} \times \mathcal{G}^{\infty}$, set

$$\theta_{\Delta}(\omega, \xi) = \begin{cases} k+1, & \text{if } \omega_0 = \Delta^k u, u \in \{a, b, ab, ba\}, \\ & u * \xi_0 \in T^{-1} \cup \Delta, \\ k, & \text{if } \omega_0 = \Delta^k u, u \in \{a, b, ab, ba\}, \\ & u * \xi_0 \notin T^{-1} \cup \Delta, \end{cases} \tag{41}$$

where $\omega = (\omega_0, \omega_1, \ldots)$. In words, $\theta_{\Delta}(\omega, \xi)$ counts the exponent of $\Delta$ created when $\xi$ is left-multiplied by $\omega_0$.

Lemma 4.10. *We have*

$$\gamma_{\Delta} = \int_{U^{\mathbb{N}} \times \mathcal{G}^{\infty}} \theta_{\Delta}(\omega, \xi)\, d\nu^{\otimes \mathbb{N}}(\omega)\, d\mu^{\infty}(\xi). \tag{42}$$

Proof. Let $\tau : (u_n)_n \mapsto (u_{n+1})_n$ be the translation shift on $U^{\mathbb{N}}$. Let $\mathcal{T}$ be the function defined by

$$\mathcal{T} : \begin{cases} U^{\mathbb{N}} \times \mathcal{G}^{\infty} \to U^{\mathbb{N}} \times \mathcal{G}^{\infty}, \\ (\omega, \xi) \to (\tau(\omega), p(\omega_0) \circledast \xi). \end{cases}$$

It follows from Proposition 4.1 that the measure $\nu^{\otimes \mathbb{N}} \otimes \mu^{\infty}$ on $U^{\mathbb{N}} \times \mathcal{G}^{\infty}$ is $\mathcal{T}$-stationary and ergodic. From the ergodic theorem, we get

$$\lim_{n \to \infty} \frac{1}{n} \sum_{i=0}^{n-1} \theta_{\Delta} \circ \mathcal{T}^i(\omega, \xi) = \int_{U^{\mathbb{N}} \times \mathcal{G}^{\infty}} \theta_{\Delta}(\omega, \xi)\, d\nu^{\otimes \mathbb{N}}(\omega)\, d\mu^{\infty}(\xi). \tag{43}$$

It remains to prove that $\gamma_{\Delta} = \lim_n n^{-1} \sum_{i=0}^{n-1} \theta_{\Delta} \circ \mathcal{T}^i(\omega, \xi)$. Define

$$X_0^{\ell} = 1, \qquad X_n^{\ell} = w_{n-1} X_{n-1}^{\ell} = w_{n-1} w_{n-2} \cdots w_0.$$

This is a realization of the *left* random walk associated with $(B_3, \nu)$, to be compared with (3). Denote the Garside normal form of $X_n^{\ell}$ by $\widehat{X}_n^{\ell} \Delta^{k_n^{\ell}}$. Obviously and almost surely, $\gamma_{\Delta} = \lim_n k_n^{\ell}/n$. [The law of $(w_{n-1}, \ldots, w_0)$ is the same as the law of $(w_0, \ldots, w_{n-1})$.] Besides, for all $m \in \mathbb{N}$, and for all $n$ large enough, the prefixes of length $m$ of $\widehat{X}_n^{\ell}$ and $w_{n-1} \circledast w_{n-2} \circledast \cdots \circledast w_0 \circledast \xi$ coincide. It follows that

$$\lim_n k_n^{\ell} - \sum_{i=0}^{n-1} \theta_{\Delta} \circ \mathcal{T}^i(\omega, \xi) \qquad \text{a.s. finite.}$$

We conclude that $\gamma_{\Delta} = \lim_n n^{-1} \sum_{i=0}^{n-1} \theta_{\Delta} \circ \mathcal{T}^i(\omega, \xi)$ and the lemma is proved. □

Using equation (42) and the relation $\mu^{\infty}(uT^{\mathbb{N}}) = R(u)$, we obtain (25).

Now let us prove (26). Recall the following notation: $(X_n)_n$ is a realization of $(B_3, \nu)$, $(Y_n)_n = (\phi \circ p(X_n))_n$ is a realization of $(B_3/Z, \mu)$, and $(Z_n)_n$ is defined as in (33). Set $Y_{\infty} = y_1 y_2 \cdots$ and $Z_{\infty} = z_1 z_2 \cdots$.

Let $\mathfrak{p}$ be the unique solution in $\{x \in (\mathbb{R}_+^*)^\Sigma | \sum_{u\in\Sigma} x(u) = 1\}$ of $xP = x$, where $P$ is defined in (22). Define $\mathfrak{p}(1) = \mathfrak{p}(a) + \mathfrak{p}(a\Delta) + \mathfrak{p}(b) + \mathfrak{p}(b\Delta)$ and $\mathfrak{p}(2) = 1 - \mathfrak{p}(1)$.

The indicator function of an event $E$ is denoted by $\mathbf{1}_E$. By the ergodic theorem for Markov chains, we have a.s.

$$\forall u \in T \qquad \lim_n \frac{1}{n}\sum_{i=1}^n \mathbf{1}_{\{z_i = u\}} = \mathfrak{p}(u). \tag{44}$$

Consider $u_1 \cdots u_k \in \mathcal{SG}$ and $w_1 \cdots w_k \in \mathcal{G}$ such that $\psi(u_1 \cdots u_k) = w_1 \cdots w_k$. Set $u_i = v_i \Delta^{0/1}, v_i \in T$. Clearly, $|v_i|_{S_+} = |w_i|_{S_+}$ for all $i$.

From this last observation and (44), we easily deduce that a.s.

$$\lim_n \frac{1}{n}\sum_{i=1}^n \mathbf{1}_{\{y_i = a \text{ or } b\}} = \mathfrak{p}(1), \qquad \lim_n \frac{1}{n}\sum_{i=1}^n \mathbf{1}_{\{y_i = ab \text{ or } ba\}} = \mathfrak{p}(2).$$

Let $Y_n[i]$ be the $i$th letter of $Y_n$. Define

$$|Y_n|_1 = \sum_{i=1}^{|Y_n|_\Sigma} \mathbf{1}_{\{Y_n[i] = a \text{ or } b\}}, \qquad |Y_n|_2 = \sum_{i=1}^{|Y_n|_\Sigma} \mathbf{1}_{\{Y_n[i] = ab \text{ or } ba\}}. \tag{45}$$

We would like to prove that a.s. $\lim_n |Y_n|_1/n = \mathfrak{p}(1)\gamma_\Sigma$ and $\lim_n |Y_n|_2/n = \mathfrak{p}(2)\gamma_\Sigma$. This requires an argument since $Y_n$ has not stabilized yet. Define

$$\forall l \in \mathbb{N} \qquad \tau_l = \inf\{N | \forall n \geq N, Y_n \in y_1 y_2 \cdots y_l \Sigma^*\}.$$

The r.v.'s $\tau_l$ are a.s. finite since the random walk is transient. Also $(\tau_l)_l$ is a strictly increasing sequence of r.v.'s and $Y_{\tau_l} = y_1 \cdots y_l$, or $Y_{\tau_l} = y_1 \cdots y_l w$, $w \in \Sigma$. Using Birkhoff's ergodic theorem, we get a.s.

$$\lim_n \frac{|Y_{\tau_n}|_1}{\tau_n} = \lim_n \frac{|Y_{\tau_n}|_1}{|Y_{\tau_n}|_\Sigma}\frac{|Y_{\tau_n}|_\Sigma}{\tau_n} = \mathfrak{p}(1)\gamma_\Sigma, \tag{46}$$

and, similarly, $\lim_n |Y_{\tau_n}|_2/\tau_n = \mathfrak{p}(2)\gamma_\Sigma$ a.s. We conclude that $\lim_n |Y_n|_1/n = \mathfrak{p}(1)\gamma_\Sigma$ and $\lim_n |Y_n|_2/n = \mathfrak{p}(2)\gamma_\Sigma$. It follows easily that (26) holds.

The only point remaining to be proved is (27). Let $\widehat{X}_n \Delta^{\delta_n}$ be the Garside normal form of $X_n$. View $\widehat{X}_n$ as a word over the alphabet $T = \{a, b, ab, ba\}$. Let $\widehat{X}_n[i]$ be the $i$th letter of $\widehat{X}_n$. Define $|\widehat{X}_n|_1$ and $|\widehat{X}_n|_2$ as in (45). (For $i = 1, 2$, we have $|\widehat{X}_n|_i = |Y_n|_i$ or $|Y_n|_i + 1$.)

Assume first that $\gamma_\Delta > 0$. Asymptotically and a.s., we have $\delta_n > 0$. Then, according to Proposition 3.2, $|X_n|_S = |\widehat{X}_n|_{S_+} + 3\delta_n$. Therefore, $\gamma = \gamma_{S_+} + 3\gamma_\Delta$.

Assume now that $\gamma_\Delta < 0, |\gamma_\Delta| < \mathfrak{p}(2)\gamma_\Sigma$. Asymptotically and a.s., we have $\delta_n < 0$ and $|\delta_n| < |\widehat{X}_n|_2$. We follow the procedure described in Proposition 3.2. When picking the largest elements in $\widehat{X}_n$, we choose only elements of length 2. When spreading the $\Delta^{-1}$'s, the transformations are consequently of the type $[ab\Delta^{-1} \to b^{-1}]$ or $[ba\Delta^{-1} \to a^{-1}]$. We deduce that $|X_n|_S = |\widehat{X}_n|_{S_+} - |\delta_n|$. Therefore, $\gamma = \gamma_{S_+} - |\gamma_\Delta|$.

Assume that $\gamma_\Delta < 0, \mathfrak{p}(2)\gamma_\Sigma < |\gamma_\Delta| < \gamma_\Sigma$. Asymptotically and a.s., we have $\delta_n < 0$ and $|\widehat{X}_n|_2 < |\delta_n| < |\widehat{X}_n|_1 + |\widehat{X}_n|_2$. When picking the largest elements in $\widehat{X}_n$, we then choose all the elements of length 2 and some of length 1. For the picked elements of length 1, the transformations are of the type $[a\Delta^{-1} \to b^{-1}a^{-1}]$ or $[b\Delta^{-1} \to a^{-1}b^{-1}]$. We deduce that $|X_n|_S = |\widehat{X}_n|_{S_+} - |\widehat{X}_n|_2 + (|\delta_n| - |\widehat{X}_n|_2)$. Hence, $\gamma = \gamma_{S_+} + |\gamma_\Delta| - 2\mathfrak{p}(2)\gamma_\Sigma$.

The last case is treated similarly. This completes the proof of (27).

4.4. *Dual structure.* Consider the following alternative presentation of the braid group on three stands:

$$\widehat{B}_3 = \langle \sigma_1, \sigma_2, \sigma_3 | \sigma_1\sigma_2 = \sigma_2\sigma_3 = \sigma_3\sigma_1 \rangle. \tag{47}$$

A possible isomorphism $\varphi : \widehat{B}_3 \to B_3$ is defined by $\varphi(\sigma_1) = a$, $\varphi(\sigma_2) = b$, $\varphi(\sigma_3) = b^{-1}ab = aba^{-1}$. There is a natural interpretation of the elements of $\widehat{B}_3$ as braids on a "cylinder," with $\sigma_3$ corresponding to the crossing of the strands 3 and 1 (cf. the Introduction).

Set $\delta = \sigma_1\sigma_2 = \sigma_2\sigma_3 = \sigma_3\sigma_1$. We have $\varphi(\delta^3) = \Delta^2$. So the center $\widehat{Z}$ of $\widehat{B}_3$ is generated by $\delta^3$. Let $\widehat{p} : \widehat{B}_3 \to \widehat{B}_3/\widehat{Z}$ be the canonical morphism. Set

$$\widehat{\Sigma} = \{\sigma_1, \sigma_2, \sigma_3, \sigma_1\delta, \sigma_2\delta, \sigma_3\delta, \sigma_1\delta^2, \sigma_2\delta^2, \sigma_3\delta^2\}. \tag{48}$$

There exists a normal form over $\widehat{\Sigma}^*$ for elements of $\widehat{B}_3/\widehat{Z}$ exactly analog to the one in Proposition 3.1.

From there, all the results in Section 4 have counterparts. So we can completely analyze the random walk $(\widehat{B}_3, \nu)$, where $\nu$ is a probability measure on $\widehat{p}^{-1}(\widehat{\Sigma} \cup \{1, \delta, \delta^2\})$ generating the whole group. Observe that

$$\widehat{\Sigma} \cup \{1, \delta, \delta^2\} \sim \{1, \Delta, a, b, ab, a\Delta, ab\Delta, ba\Delta, a \cdot ab, a \cdot a\Delta, ab \cdot b, ab \cdot ba\Delta\}.$$

In particular, the sets $\varphi(\widehat{\Sigma} \cup \{1, \delta, \delta^2\})$ and $\Sigma \cup \{1, \Delta\}$ (viewed as subsets of $B_3/Z$) are incomparable for inclusion. So the random walks that we are able to treat by working with $B_3$ or $\widehat{B}_3$ are indeed different.

The presentation in (47) is an instance of the so-called *dual structure* for Artin groups of finite Coxeter type. The dual structure is linked to Garside monoids and groups. This is a subject of recent but intense research; see, for instance, [5, 28].

**5. Artin groups of dihedral type.** All the above results adapt to any Artin group of dihedral type. Define $\mathrm{prod}(u, v; k) = uvuv \cdots$, with $k$ letters in the word on the right-hand side. The *Artin groups of dihedral type* are defined by the group presentations, for $k \geq 3$,

$$A_k = \langle a, b | \mathrm{prod}(a, b; k) = \mathrm{prod}(b, a; k) \rangle.$$

Observe that $A_3 = B_3$. Now consider $A_k$ for a fixed $k$. Set $\Delta = \mathrm{prod}(a, b; k) = \mathrm{prod}(b, a; k)$. If $k$ is even, then the center of $A_k$ is $Z = \langle \Delta \rangle$; if $k$ is odd, then the center of $A_k$ is $Z = \langle \Delta^2 \rangle$. Consider $A_k/Z$. Define the following subsets of $A_k/Z$:

$$S = \{a, b, a^{-1}, b^{-1}\},$$
$$T = \{\mathrm{prod}(a, b; i), \mathrm{prod}(b, a; i);\ i = 1, \dots, k-1\}, \qquad \Sigma = T \cup T^{-1}.$$

If $k$ is even, then $\Sigma = T$, and if $k$ is odd, then $\Sigma = T \sqcup T^{-1} = T \sqcup \{u\Delta, u \in T\}$. The exact analog of Proposition 3.1 holds for $A_k$ and $A_k/Z$. In particular, for $g \in A_k/Z$, there exist $g_1, \dots, g_m \in T$ such that

$$g = \begin{cases} g_1 \cdots g_m, & \text{if } k \text{ is even,} \\ g_1 \cdots g_m \Delta^n, n \in \{0, 1\}, & \text{if } k \text{ is odd,} \end{cases}$$

and the decomposition is unique if $m$ is chosen to be minimal. Like in Section 3, this defines the set $\mathcal{G} \subset T^* \Sigma$ of Garside normal forms of elements of $A_k/Z$. The set $\mathcal{G}$ forms a geodesic cross-section of $A_k/Z$ over $\Sigma$. Denote by $\mathcal{SG} \subset \Sigma^*$ the set of all the geodesics. If $k$ is even, then $\mathcal{G} = \mathcal{SG}$ and if $k$ is odd, then $\mathcal{G} \subsetneq \mathcal{SG}$. Define the maps $\mathrm{First}, \mathrm{Last} : \Sigma \to \{a, b\}$ as follows:

For $u \in \Sigma$, $\mathrm{First}(u)$ is the first symbol of $u$; for $k$ even and $u \in \Sigma$, $\mathrm{Last}(u)$ is the last symbol of $u$; for $k$ odd and $u \in T$, $\mathrm{Last}(u)$ is the last symbol of $u$; for $k$ odd and $u\Delta \in T^{-1}$, $\mathrm{Last}(u\Delta)$ is the last symbol of $v$ where $\Delta v = u\Delta$.

The characterizations (7) and (12) and the properties (8) and (11) still hold [for $k$ even, keep only the lines 1, 3 and 5 in (11)]. For $k$ even, the map $\psi : \mathcal{SG} \to \mathcal{G}$ is just the identity. For $k$ odd, the map $\psi : \mathcal{SG} \to \mathcal{G}$ is still recognizable by a transducer, and the transducer recognizing it is the natural generalization of the one in Figure 8. The natural generalization of Proposition 3.2 holds, see Proposition 5.3 below and [22] for a proof.

So all the ingredients are in place to proceed as in Section 4. Set $\overline{\Sigma} = \Sigma \cup \{1\}$ if $k$ is even and $\overline{\Sigma} = \Sigma \cup \{1, \Delta\}$ if $k$ is odd. Let $\mu$ be a probability measure on $\overline{\Sigma}$ such that $\bigcup_n \mathrm{supp}(\mu^{*n}) = A_k/Z$. The end boundary of $\mathcal{X}(A_k/Z, \Sigma)$ is still $\mathcal{G}^\infty$ defined by (15). Let $\mu^\infty$, defined on $\mathcal{G}^\infty$, be the harmonic measure of $(A_k/Z, \mu)$. For $u \in \Sigma$, define $\mathrm{Next}(u) = \{v \in \Sigma \,|\, \mathrm{Last}(u) = \mathrm{First}(v)\}$. The *traffic equations* of $(A_k/Z, \mu)$ are defined by (19) with $[v, w \in \overline{\Sigma}]$ instead of $[v, w \in \Sigma]$, with $\iota$ defined as in (10) if $k$ is odd, and with $[\iota(u)\Delta]$ replaced by $[u]$ if $k$ is even.

Let $p : A_k \to A_k/Z$ be the canonical morphism. We have the following result:

Proposition 5.1. *Consider the Artin group $A_k = \langle a, b \,|\, \mathrm{prod}(a, b; k) = \mathrm{prod}(b, a; k) \rangle$. Consider the random walk $(A_k, \nu)$, where $\nu$ is a probability measure on $p^{-1}(\overline{\Sigma}) = \{u\Delta^n, u \in T \cup \{1\}, n \in \mathbb{Z}\}$ such that $\bigcup_n \mathrm{supp}(\nu^{*n}) = A_k$ and $\sum_x |x|_S \nu(x) < \infty$. Set $\mu = \nu \circ p^{-1}$. The traffic equations of $(A_k/Z, \mu)$*

*have a unique solution* $r \in \{x \in (\mathbb{R}_+^*)^\Sigma \mid \sum_{u\in\Sigma} x(u) = 1\}$. *Set* $q(u) = r(u)/r(\mathrm{Next}(u))$ *for* $u \in \Sigma$. *The harmonic measure* $\mu^\infty$ *of* $(A_k/Z, \mu)$ *is given by,* $\forall u = u_1 \cdots u_n \in \mathcal{G} \cap T^n$,

$$\mu^\infty(uT^{\mathbb{N}}) = \begin{cases} q(u_1)\cdots q(u_{n-1})r(u_n), & \text{if } k \text{ is even,} \\ \displaystyle\sum_{v_1\dots v_n \in \psi^{-1}(u)} q(v_1)\cdots q(v_{n-1})[r(v_n) + r(v_n\Delta)], & \text{if } k \text{ is odd.} \end{cases}$$

For $u \in T$, set $R(u) = \mu^\infty(uT^{\mathbb{N}})$. For $k$ even, we have $R(u) = r(u)$, and for $k$ odd, we have $R(u) = r(u) + r(u\Delta)$. Let $\ell : \Sigma \to T$ be defined as follows. If $k$ is odd, set $\ell(u) = u$ if $u \in T$, and $\ell(u) = v$ if $u = v\Delta \in T^{-1}$. If $k$ is even, set $\ell(u) = u$. Define $P \in \mathbb{R}_+^{\Sigma\times\Sigma}$ as in (22). Let $\mathfrak{p} \in \{x \in (\mathbb{R}_+^*)^\Sigma \mid \sum_{u\in\Sigma} x(u) = 1\}$ be the unique solution to the equations $xP = x$. Set $S^+ = \{a, b\}$. For $i \in \{1, \dots, k-1\}$, set

$$\text{(49)}\qquad \begin{aligned} &\text{for } k \text{ even,} \qquad \mathfrak{p}(i) = \sum_{u\in T, |u|_{S_+} = i} \mathfrak{p}(u), \\ &\text{for } k \text{ odd,} \qquad \mathfrak{p}(i) = \sum_{u\in T, |u|_{S_+} = i} \mathfrak{p}(u) + \mathfrak{p}(u\Delta). \end{aligned}$$

Use the convention $\sum_{j=1}^{0} * = 0$. Define the drifts $\gamma$, $\gamma_\Sigma$, $\gamma_{S_+}$ and $\gamma_\Delta$, as in (14). We have:

Proposition 5.2. *For $k$ odd, the drifts $\gamma_\Sigma$, $\gamma_{S_+}$ and $\gamma_\Delta$ are given respectively by (23), (24) and (25). For $k$ even, the drifts $\gamma_\Sigma$ and $\gamma_{S_+}$ are given by (23) and (24) under the convention that $\forall u, \mu(\Delta u^{-1}) = 0$. For $k$ even, the drift $\gamma_\Delta$ is given by $\gamma_\Delta = \sum_{u\in T\cup\{1\}, n\in\mathbb{Z}} n\nu(u\Delta^n)$. The following relation holds for all $k$: $\gamma_{S_+} = \sum_{i=1}^{k-1} i\mathfrak{p}(i)\gamma_\Sigma$. For all $k$, the drift $\gamma$ is given by*

$$\text{(50)}\qquad \gamma = \begin{cases} \gamma_{S_+} + k\gamma_\Delta, & \text{if } \gamma_\Delta \geq 0, \\ \gamma_{S_+} - (k-2i)|\gamma_\Delta| - \displaystyle\sum_{j=1}^{i-1} (2i-2j)\mathfrak{p}(k-j)\gamma_\Sigma, & \text{if } \gamma_\Delta < 0, \end{cases}$$

*for $i = 1 + \max\{\ell \in \{0, \dots, k-1\} \mid \sum_{j=1}^{\ell} \mathfrak{p}(k-j)\gamma_\Sigma < |\gamma_\Delta|\}$. Observe that if $\gamma_\Delta < 0, \gamma_\Sigma < |\gamma_\Delta|$, then (50) reduces to $\gamma = -\gamma_{S_+} + k|\gamma_\Delta|$.*

For $k$ even, we have $(Z_n)_n = (Y_n)_n$ since $\mathcal{G} = \mathcal{SG}$. So the approach followed in Section 4 degenerates and coincides precisely with the simpler approach used for so-called 0-automatic pairs in [19]. Indeed, for $k$ even, the pair $(A_k/Z, \Sigma)$ is 0-automatic. For $k$ odd, the pair $(A_k/Z, \Sigma)$ is not 0-automatic and introducing the process $(Z_n)_n \neq (Y_n)_n$ is the only way to proceed. We now detail the proof.

Proof of Proposition 5.2. The formulae (23), (24) and (25) for the drifts $\gamma_\Sigma$, $\gamma_{S_+}$ and $\gamma_\Delta$ are proven exactly as in Section 4 (with a simplification if $k$ is even due to the fact that $\mathcal{G} = \mathcal{SG}$). Only equation (50) requires some details.

We first explain how to compute the length of a group element $g \in A_k$ with the help of its Garside normal form.

Define the following subset of $A_k$:

$$T = \{\mathrm{prod}(a, b; i), \mathrm{prod}(b, a; i);\ \ i = 1, \ldots, k-1\}.$$

The Garside normal form of $g \in A_k$ is the unique decomposition $g = g_1 \cdots g_m \Delta^\delta$ with $g_1, \ldots, g_m \in T$, $\delta \in \mathbb{Z}$ and $m$ minimal. Write

$$\overline{g} = \Delta g \Delta^{-1} = \Delta^{-1} g \Delta$$

(since $\Delta^2$ is central) for the conjugate of $g$ by $\Delta$. For instance, if $k$ is odd, one has $\overline{a} = b$ and $\overline{b} = a$. If $k$ is even, $\Delta$ is central so $\overline{a} = a$ and $\overline{b} = b$. The conjugacy by $\Delta$ is an involution of $A_k$, which leaves invariant both the set $T$ and the length $|\cdot|_S$. Note that we have

$$\forall n \in \mathbb{N}^* \qquad g\Delta^n = \begin{cases} \Delta^n \overline{g}, & \text{if } n \text{ is odd,} \\ \Delta^n g, & \text{if } n \text{ is even.} \end{cases} \tag{51}$$

Last, for $g \in T$, we set $g^* = \Delta g^{-1}$. Note that $g^* \in T$, that $\overline{g^*} = \overline{g}^* = \overline{g}^{-1}\Delta$. Moreover, for every $g \in T$,

$$\mathrm{Last}(g^*) = \begin{cases} a, & \text{if } \mathrm{First}(g) = b, \\ b, & \text{if } \mathrm{First}(g) = a. \end{cases} \tag{52}$$

Keeping in mind that $g^* \cdot g = \Delta$, one may think of $g^*$ as the left-complement of $g$ to get $\Delta$. For instance, in the braid group $B_3$, we have $\Delta = aba = bab$, therefore, $a^* = ab$ and $(ab)^* = b$.

The next observation is central in what follows:

$$\forall g \in T \qquad |g^*|_S = |\Delta|_S - |g|_S = k - |g|_S. \tag{53}$$

The Proposition 5.3 below is a restatement of Proposition 3.2 in the case of Artin groups $A_k$. See [22] for a proof. It can be viewed as a "Greedy Algorithm" to get a minimal length representative starting from a Garside normal form.

Proposition 5.3. *Let $g_1 \cdot \cdots \cdot g_m \cdot \Delta^\delta$ be the Garside normal form of an element $g$ in $A_k$.*

- *Assume $\delta \geq 0$. Then the word $g_1 \cdot \cdots \cdot g_m \cdot \Delta^\delta$ is of minimal length with respect to $S$, and we have*

$$|g|_S = \sum_{i=1}^{m} |g_i|_S + \delta \cdot k. \tag{54}$$

- *Assume $\delta \leq -m$ and write, for $1 \leq j \leq m$, $h_j = (g_j^*)^{-1}$ if $j$ is odd, and $h_j = (\overline{g}_j^*)^{-1}$ if $j$ is even. Then the word $h_1 \cdot \cdots \cdot h_m \cdot \Delta^{-(|\delta|-m)}$ is a word of minimal length representing $g$, and we have*

$$|g|_S = \sum_{i=1}^{m} |h_i|_S + |\Delta^{-(|\delta|-m)}|_S = \sum_{i=1}^{m} (k - |g_i|_S) + k(|\delta| - m). \tag{55}$$

- *Assume $-m < \delta < 0$, and choose among $g_1, \ldots, g_m$ a family $(g_i)_{i \in L}$ of $|\delta|$ elements of maximal possible length. We define $h_1, \ldots, h_m \in T \sqcup T^{-1}$ by*

$$h_j = \begin{cases} g_j, & \textit{if } j \notin L \textit{ and } \#(\{1, \ldots, j\} \cap L) \textit{ is even,} \\ \overline{g}_j, & \textit{if } j \notin L \textit{ and } \#(\{1, \ldots, j\} \cap L) \textit{ is odd,} \\ (\overline{g}_j^*)^{-1}, & \textit{if } j \in L \textit{ and } \#(\{1, \ldots, j\} \cap L) \textit{ is even,} \\ (g_j^*)^{-1}, & \textit{if } j \in L \textit{ and } \#(\{1, \ldots, j\} \cap L) \textit{ is odd.} \end{cases} \tag{56}$$

*Then the word $h_1 \cdot \cdots \cdot h_m$ is a word of minimal length representing $g$, and we have*

$$|g|_S = \sum_{i \in L} (k - |g_i|_S) + \sum_{i \notin L} |g_i|_S. \tag{57}$$

Let us now prove equation (50). Denote by $\widehat{X}_n \Delta^{\delta_n}$ the Garside normal form of $X_n$. View $\widehat{X}_n$ as a word over the alphabet $T$. Let $\widehat{X}_n[i]$ be the $i$th letter of $\widehat{X}_n$ and set, for $j \in \{1, \ldots, k-1\}$,

$$|\widehat{X}_n|_j = \sum_{i=1}^{|\widehat{X}_n|_T} \mathbf{1}_{\{\widehat{X}_n[i] = \mathrm{prod}(a,b;j) \text{ or } \mathrm{prod}(b,a;j)\}}.$$

Obviously, $|\widehat{X}_n|_T = \sum_{j=1}^{k-1} |\widehat{X}_n|_j$. Consider $\mathfrak{p}$ defined in equation (49). Observe that $\mathfrak{p}(j)$ is a.s. the asymptotic proportion of letters of length $j$ occurring in $\widehat{X}_n$. Moreover, we have (recall the proof of Theorem 4.3)

$$\forall j \in \{1, \ldots, k-1\} \qquad \lim_n \frac{|\widehat{X}_n|_j}{n} = \mathfrak{p}(j)\gamma_\Sigma. \tag{58}$$

Assume first that $\gamma_\Delta > 0$. Asymptotically and a.s., we have $\delta_n > 0$. Then, according to equation (54), $|X_n|_S = |\widehat{X}_n|_{S_+} + k\delta_n$. Therefore, $\gamma = \gamma_{S_+} + k\gamma_\Delta$.

Assume now that $\gamma_\Delta < 0, |\gamma_\Delta| < \gamma_\Sigma$. Let $i \in \{1, \ldots, k-1\}$ the (unique) integer satisfying

$$\sum_{j=1}^{i-1} \mathfrak{p}(k-j)\gamma_\Sigma < |\gamma_\Delta| < \sum_{j=1}^{i} \mathfrak{p}(k-j)\gamma_\Sigma.$$

Asymptotically and almost surely, we have $\delta_n < 0$ and

$$\sum_{j=0}^{i-1} |\widehat{X}_n|_{k-j} < |\delta_n| < \sum_{j=0}^{i} |\widehat{X}_n|_{k-j}.$$

Now, we follow the procedure described in Proposition 5.3. When picking the $|\delta_n|$ largest elements in $\widehat{X}_n$, we choose all the elements of length $k-1, k-2, \dots, k-i+1$ and some of length $k-i$. According to equation (57), we have

$$|X_n|_S = \sum_{j=1}^{i-1} j\,|\widehat{X}_n|_{k-j} + \sum_{j=i+1}^{k-1} (k-j)|\widehat{X}_n|_{k-j} \\ + i\left(|\delta_n| - \sum_{j=1}^{i-1} |\widehat{X}_n|_{k-j}\right) + (k-i)\left(\sum_{j=1}^{i} |\widehat{X}_n|_{k-j} - |\delta_n|\right).$$

Rearranging this sum and using the fact that $|X_n|_{S^+} = \sum_{j=1}^{k-1}(k-j)|\widehat{X}_n|_{k-j}$, we get

$$|X_n|_S = |X_n|_{S^+} + (2i-k)|\delta_n| + \sum_{j=1}^{k-1}(2j-2i)|\widehat{X}_n|_{k-j}.$$

Therefore, we have $\gamma = \gamma_{S_+} - (k-2i)|\gamma_\Delta| - \sum_{j=1}^{i-1}(2i-2j)\mathfrak{p}(k-j)\gamma_\Sigma$, which is equation (50).

Last, assume that $\gamma_\Delta < 0, |\gamma_\Delta| > \gamma_\Sigma$. Asymptotically and a.s., we have $|\widehat{X}_n|_T < |\delta_n|$. Then, equation (55) gives

$$|X_n|_S = \sum_{j=1}^{k-1}(k-j)|\widehat{X}_n|_j + k\,(|\delta_n| - |\widehat{X}_n|_T) = -|\widehat{X}_n|_{S^+} + k|\delta_n|,$$

and, consequently, $\gamma = -\gamma_{S_+} + k|\gamma_\Delta|$. □

*Dual structure.* The dual structure of $A_k$ [cf. Section 4.4] is given by (cf., e.g., [28], Chapter V)

$$\widehat{A}_k = \langle \sigma_1, \dots, \sigma_k | \sigma_1\sigma_2 = \sigma_2\sigma_3 = \cdots = \sigma_k\sigma_1 \rangle.$$

Set $\delta = \sigma_1\sigma_2$. There exist counterparts of Proposition 5.1 and 5.2 with $\widehat{A}_k$ instead of $A_k$ and $\widehat{\Sigma} = \{\sigma_i\delta^j, \delta^j, i \in \{1,\dots,k\}, j \in \{0,\dots,k-1\}\}$ instead of $\overline{\Sigma}$.

5.1. *Simple random walk.* Let us consider the simple random walk $(A_k, \nu)$ defined by $\nu(a) = \nu(a^{-1}) = \nu(b) = \nu(b^{-1}) = 1/4$.

Let $r$ be the unique solution to the traffic equations. By obvious symmetry, we must have $r(\mathrm{prod}(a,b,i)) = r(\mathrm{prod}(b,a,i))$ for all $i \in \{1,\dots,k-1\}$ if $k$ is even, and $r(u) = r(\iota(u)\Delta)$ for all $u \in T$ if $k$ is odd. Simplifying the traffic equations accordingly gives precisely the traffic equations for the simple random walk on $\mathbb{Z}/k\mathbb{Z} \star \mathbb{Z}/k\mathbb{Z}$ with respect to a minimal set of generators. This reduction can also be understood geometrically on the Cayley graphs as illustrated in Figure 10. The explicit computations for $\mathbb{Z}/k\mathbb{Z} \star \mathbb{Z}/k\mathbb{Z}$ are carried out in [21], Section 4.4. This leads to the results below (details can be found in [23], Section 4).

Consider the applications $F_n:[0,1]\to\mathbb{R}, n\in\mathbb{N}$, defined by

$$F_0(x)=1, \qquad F_1(x)=x, \qquad \forall n\geq 2 \qquad (59)$$
$$F_n(x)=2(2-x)F_{n-1}(x)-F_{n-2}(x).$$

For $k\geq 3$, the equation $F_k(x)=1$ has a unique solution in $(0,1)$ that we denote by $x_k$; see [21], Lemma 4.1.

Proposition 5.4. *Consider the simple random walk $(A_k,\nu)$ with $\nu(a)=\nu(a^{-1})=\nu(b)=\nu(b^{-1})=1/4$. The drifts $\gamma_\Sigma$ and $\gamma_\Delta$ are given by*

$$\gamma_\Sigma=\frac{1-x_k}{2}, \qquad \gamma_\Delta=-\frac{1-x_k}{4}. \qquad (60)$$

*Let $\gamma$ be the drift of the length with respect to the natural generators $\{a,b,a^{-1}, b^{-1}\}$. We have*

$$\gamma=\begin{cases}(1-x_k)\left[\sum_{i=1}^{j-1} iF_i(x_k)+(j/2)F_j(x_k)\right], & \text{if } k=2j,\\ (1-x_k)\left[\sum_{i=1}^{j} iF_i(x_k)\right], & \text{if } k=2j+1.\end{cases} \qquad (61)$$

Both $\gamma_\Sigma$ and $\gamma$ are increasing functions of $k$, while $\gamma_\Delta$ is a decreasing function of $k$. We have $\lim_k\gamma_\Sigma=1/3$, $\lim_k\gamma_\Delta=-1/6$ and $\lim_k\gamma=1/2$.

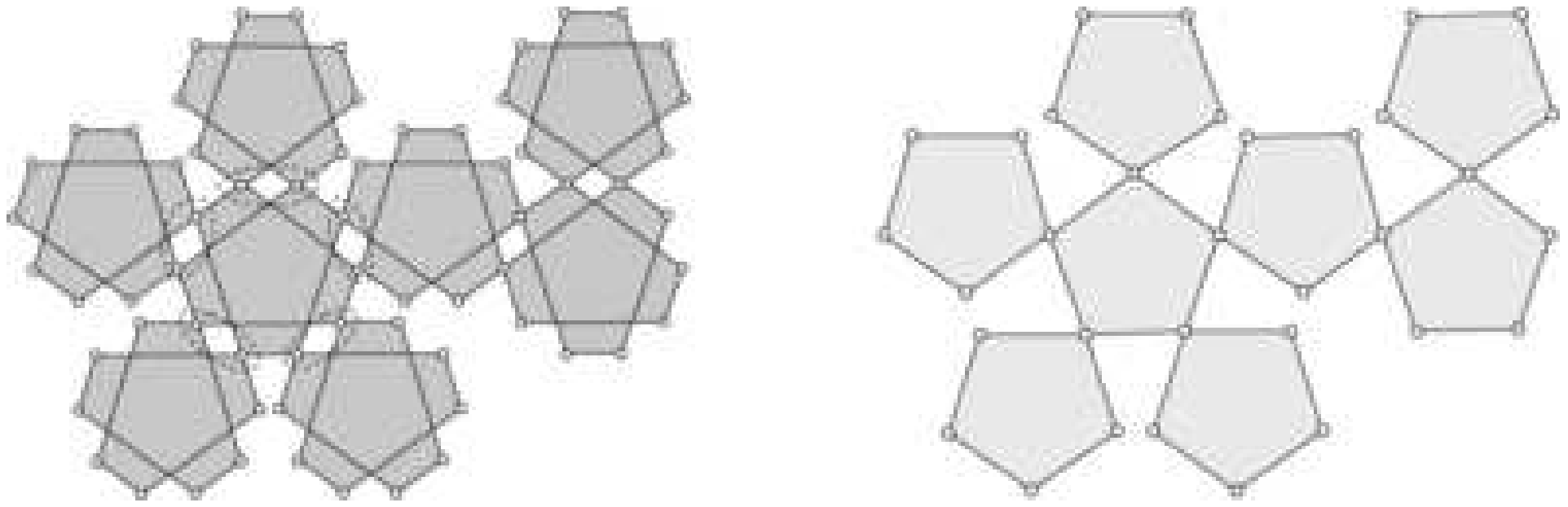

Fig. 10. *The Cayley graphs $\mathcal{X}(A_5/Z,S)$ (left) and $\mathcal{X}(\mathbb{Z}/5\mathbb{Z}\star\mathbb{Z}/5\mathbb{Z},\widetilde{S})$ (right).*

Table 2

| | $A_3$ | $A_4$ | $A_5$ | $A_6$ | $A_7$ | $A_8$ |
|---|---|---|---|---|---|---|
| $\gamma$ | $1/4$ | $(\sqrt{5}-1)^2/4$ | $(\sqrt{13}-1)^2/16$ | 0.462598... | 0.475221... | 0.487636... |

In Table 2 there are the first values of $\gamma$, given either in closed form or numerically when no closed form could be found.

**Acknowledgments.** The authors are grateful to Christian Blanchet, Yves Guivarc'h, Matthieu Picantin, Alain Valette and Dimitri Zvonkine for helpful discussions. Special thanks go to Neil O'Connell for suggesting the manta ray likeness.

LIAFA
CNRS–Université Paris 7
case 7014
2 place Jussieu
75251 Paris Cedex 05
France
E-mail: Jean.Mairesse@liafa.jussieu.fr
URL: http://www.liafa.jussieu.fr/~mairesse/

LMAM
Université de Bretagne-Sud
Campus de Tohannic
BP 573
56017 Vannes
France
E-mail: Frederic.Matheus@univ-ubs.fr
URL: http://web.univ-ubs.fr/lmam/matheus/